\documentclass[letterpaper, 12 pt]{article}
\usepackage{amsthm, amsmath, amssymb, amsfonts, graphicx}

\def\EE{\mathcal{E}}
\def\FF{\mathcal{F}}
\def\LL{\mathcal{L}}
\def\II{\mathcal{I}}
\def\UU{\mathcal{U}}
\def\GG{\mathcal{G}}

\def\Copen{\mathring{C}}
\def\Cbound{\partial C}

\def\ZZ{\mathbb{Z}}
\def\RR{\mathbb{R}}

\def\C{\mathbf{C}}
\def\N{\mathbf{N}}
\def\E{\mathbf{E}}
\def\S{\mathbf{S}}
\def\W{\mathbf{W}}
\def\NE{\N\E}
\def\NW{\N\W}
\def\NW{\N\W}
\def\SE{\S\E}
\def\SW{\S\W}
\def\q{\mathbf{q}}

\def\cont{\ensuremath{\Rightarrow \Leftarrow}}

\def\maps{\colon}
\def\into{\hookrightarrow}

\def\gwof{graph with open faces }
\def\gwofs{graphs with open faces }

\def\Mout{M_{\textrm{out}}}
\def\Gout{G_{\textrm{out}}}

\def\cross{\includegraphics[width=0.2 in]{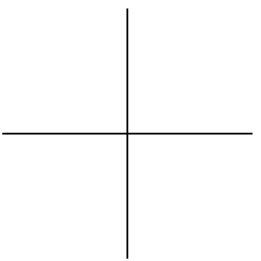}}

\def\backwrench{\includegraphics[width=0.2 in]{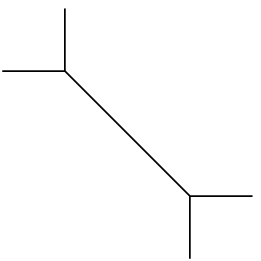}}

\def\wrench{\includegraphics[width=0.2 in]{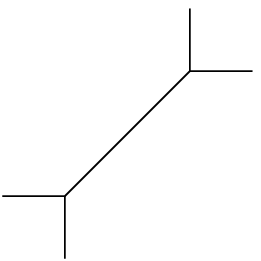}}

\newtheorem*{AD}{The Aztec Diamond Theorem}
\newtheorem*{MT1}{The Main Theorem (First Statement)}
\newtheorem*{MT2}{The Main Theorem (Second Statement)}
\newtheorem*{MT}{The Main Theorem}
\newtheorem*{MTInfExt}{The Main Theorem (Using Infinite Completions)}

\newtheorem*{URT}{The Urban Renewal Theorem}
\newtheorem*{Kuo}{Kuo's Condensation Theorem}
\newtheorem*{Heights}{The Correspondence Between Heights and Matchings}

\newtheorem{Definition}{Definition}
\newtheorem{Proposition}[Definition]{Proposition}
\newtheorem{Lemma}[Definition]{Proposition}
\newtheorem{Case}{Case}
\newtheorem{cor}[Definition]{Corollary}

\newcommand{\setleftmargin}[1]{
        \addtolength{\textwidth}{\oddsidemargin}
        \addtolength{\textwidth}{1in}
        \addtolength{\textwidth}{-#1}
        \setlength{\oddsidemargin}{-1in}
        \addtolength{\oddsidemargin}{#1}
        \setlength{\evensidemargin}{\oddsidemargin}
}
\newcommand{\setrightmargin}[1]{
        \setlength{\textwidth}{8.5in}
        \addtolength{\textwidth}{-\oddsidemargin}
        \addtolength{\textwidth}{-1in}
        \addtolength{\textwidth}{-#1}
}

\setleftmargin{1 in}
\setrightmargin{1 in}

\begin{document}

\title{Perfect Matchings and The Octahedron Recurrence}
\author{David E Speyer}

\maketitle

\begin{abstract}
We study a recurrence defined on a three dimensional lattice and prove that its values are Laurent polynomials in the initial conditions with all coefficients equal to one. This recurrence was studied by Propp and by Fomin and Zelivinsky. Fomin and Zelivinsky were able to prove Laurentness and conjectured that the coefficients were 1. Our proof establishes a bijection between the terms of the Laurent polynomial and the perfect matchings of certain graphs, generalizing the theory of Aztec diamonds. In particular, this shows that the coefficients of this polynomial, and polynomials obtained by specializing its variables, are positive, a conjecture of Fomin and Zelevinsky.
\end{abstract}

\section{Historical Introduction}

The octahedron recurrence is the product of three chains of research seeking a common generalization. The first is the study of the algebraic relations between the connected minors of a matrix, and particularly of a recurrence relating them known and Dodgson condensation. Attempting to understand the combinatorics of Dodgson condensation lead to the discovery of alternating sign matrices and Aztec diamonds. Aztec diamonds are graphs whose perfect matchings have extremely structured combinatorics and soon formed their own, second line of research as other graphs were discovered with the similar regularities. The third is the study of Somos-sequences and the Laurent phenomenon, which began with an attempt to understand theta functions from a combinatorial perspective. In this section we will sketch both lines of research. In the following section, we will begin to describe the vocabulary and main results of this paper.

\subsection{Dodgson Condensation and the Octahedron Recurrence}

Let $m_{ij}$ be a collection of formal variables indexed by ordered pairs $(i,j) \in \ZZ^2$. Let $i_1$, $i_2$, $j_1$ and $j_2$ be integers with $i_2-i_1=j_2-j_1>0$ and set
$$D_{i_1, j_1}^{i_2, j_2} = \det \left( m_{i,j} \right)_{\begin{substack} i_1 \leq i \leq i_2, \\ j_1 \leq j \leq j_2 \end{substack}}.$$
Charles Dodgson~\cite{Dodg} observed that
$$D_{i_1, j_1}^{i_2, j_2} D_{(i_1 -1), (j_1-1)}^{(i_2+1),(j_2+1)} = D_{i_1, j_1}^{(i_2+1), (j_2+1)} D_{(i_1 -1), (j_1-1)}^{i_2, j_2} - D_{(i_1-1), j_1}^{i_2, (j_2+1)} D_{i_1, (j_1-1)}^{(i_2+1), j_2}$$
suggested using this recursion to compute determinants and found generalizations that could be used when the above recursion requires dividing by $0$. (This formula is not original to Dodgson, according to page 7 of \cite{Lun} this formula has been variously attributed to Sylvester, Jacobi, Desnanot, Dodgson and Frobenius. However, the use of this identity iteratively to compute connected minors and the detailed study of this method appears to have begun with Dodgson and is known as Dodgson condensation.)

Before passing on, we should note that the study of algebraic relations between determinants and how they are effected by various vanishing conditions is essentially the study of the flag manifold and Schubert varieties. This study lead to the invention of cluster algebras which later returned in the Laurentness proofs we will discuss in Section \ref{Somos}. See \cite{FomZel2} for the definitions of Cluster Algebras. The precise relation between cluster algebras and (double) Schubert varieties is discussed in \cite{BFZ}, although the reader may wish to first consult the references in that paper for the history of results that preceeded it.

One can visualize this recurrence as taking place on a three dimensional lattice. The values $m_{ij}$ are written on a horizontal two dimensional plane. The values $D_{i_1 j_1}^{i_2 j_2}$ live on a larger three dimensional lattice with $D_{i_1 j_1}^{i_2 j_2}$ written above the center of the matrix of which $D_{i_1 j_1}^{i_2 j_2}$ is the determinant and at a height proportional to $i_2-i_1=j_2-j_1$. The six entries involved in computing the recurrence lie at the vertices of an octahedron. We make a change of the indexing variables of the lattice so that the lattice consists of those $(n,i,j)$ with $n+i+j \equiv 0 \mod 2$. Then the recurrence above can be rewritten as
$$f(n,i,j) f(n-2,i,j) = f(n-1,i-1,j) f(n-1,i+1,j) - f(n-1,i,j-1) f(n-1,i,j+1)$$
with $f(0,i,j)=m_{\frac{i+j}{2} \frac{i-j}{2}}$ and $f(-1,i,j)=1$.

Robbins and Rumsey \cite{RobRum} considered modifying the recurrence to 
$$f(n,i,j) f(n-2,i,j) = f(n-1,i-1,j) f(n-1,i+1,j) + \lambda f(n-1,i,j-1) f(n-1,i,j+1)$$
and also generalized by allowing $f(-1,i,j)$ to be arbitrary instead of being forced to be $1$. We write $x(i,j)=f(0,i,j)$ when $i+j$ is even and $x(i,j)=f(-1,i,j)$ when $i+j$ is odd. Robbins and Rumsey discovered and proved that, despite these modifications, the functions $f(n,i,j)$ were still Laurent polynomials in the $x(i,j)$ and in $\lambda$, with all coefficients equal to $1$.

Robbins and Rumsey discovered that the exponents of the $x(i,j)$ in any monomial occuring in a $f(n_0,i_0,j_0)$ formed pairs of ``compatible alternating sign matrices.'' We will omit the definition of these obects in favor of discussing a simpler combinatorial description, found by Elkies, Kuperberg, Larsen and Propp \cite{EKLP}, in terms of tilings of Aztec diamonds.

\subsection{Aztec Diamonds} \label{AztecIntro}

If $G$ is a graph then a \emph{matching} of $G$ is defined to be a collection $M$ of edges of $G$ such that each vertex of $G$ lies on exactly one edge of $M$. Matchings are also sometimes known as ``perfect matchings'', ``dimer covers'' and ``1-factors''.  Fix integers $n_0$, $i_0$ and $j_0$ with $n_0+i_0+j_0 \cong 0 \mod 2$. In this section, we will describe a combinatorial interpretation of the entries of $f(n_0,i_0,j_0)$ in terms of matchings of certain graphs called Aztec diamonds. For simplicity, we take $\lambda=1$ in the recurrence for $f$.

Consider the infinite square grid graph: its vertices are at the points of the form $(i+1/2, j+1/2)$ with $(i,j) \in \ZZ^2$, its edges join those vertices that differ by $(0, \pm 1)$ or $(\pm 1, 0)$ and its faces are unit squares centered at each $(i,j) \in \ZZ^2$. We will refer to the face centered at $(i,j)$ as $(i,j)$. Consider the set $S(n_0,i_0,j_0)$ of all faces $(i,j)$ with $|i-i_0|+|j-j_0|<n_0$. Let $G(n_0,i_0,j_0)$ be the induced subgraph of the square grid graph whose vertices are adjacent to some face of $S(n_0,i_0,j_0)$. We will call this the Aztec diamond of order $n_0$ centered at $(i_0,j_0)$.

We define the \emph{faces} of $G(n_0,i_0,j_0)$ to be the squares $(i,j)$ for $|i-i_0|+|j-j_0| \leq n_0$. Note that this includes more squares than those in $S(n_0,i_0,j_0)$; in particular, it contains some squares only part of whose boundary lies in $G(n_0,i_0,j_0)$. We refer to the faces in $S(n_0,i_0,j_0)$ as \emph{closed} faces and the others as \emph{open} faces. Let $M$ be a matching of $G(n_0,i_0,j_0)$. If $(i,j)$ is a closed face then either $0$, $1$ or $2$ edges of $(i,j)$ are used in $M$; we define $\epsilon(i,j)$ to be $-1$, $0$ or $1$ respectively. If $x(i,j)$ is an open face then either $0$ or $1$ edges will be used and, we define $\epsilon(i,j)$ to be $0$ or $1$ respectively. Set $m(M)=\prod x(i,j)^{\epsilon(i,j)}$, where the product is over all faces $(i,j)$ of $G(n_0,i_0,j_0)$.

\begin{AD}
$f(n_0,i_0,j_0)=\sum m(M)$, where the sum is over all matchings $M$ of $G(n_0,i_0,j_0)$.
\end{AD}

\begin{proof}
It seems to be difficult to find a reference which states this result exactly in this form. \cite{RobRum} shows that $f(n_0,i_0,j_0)$ is equal to a sum over ``compatible pairs of alternating sign matrices'' and \cite{EKLP} shows that such pairs are in bijection with perfect matchings of the Aztec Diamond. Tracing through these bijections easily gives the claimed theorem. Another proof can be found by mimicing the method used in \cite{Kuo} to count the matchings of the Aztec Diamond. This result will also be a corollary of this paper's main result. 
\end{proof}

\begin{cor}
The order $n$ Aztec diamond has $2^{n(n+1)/2}$ perfect matchings. 
\end{cor}

\begin{proof}
Clearly, if we plug in $1$ for each $x(i,j)$, then $m(M)=1$ for every matching $M$ so the number of matchings of the order $n$ Aztec diamond will be $f(n,i,j)|_{x(i,j)=1}$. We will denote this number by $g(n)$. Then the octahedron recurrence gives $g(n) g(n-2)=2 g(n-1)^2$ and the claim follows by induction.
\end{proof}

Many other families of graphs have been discovered for which the number of matchings of the $n^{\textrm{th}}$ graph is proportional to a constant raised to a quadratic in $n$. They are summarized and unified in \cite{Ciucu}. One of the purposes of this paper is to give proofs of these formulae as simple as the corrolary above.

The monomial $m(M)$ keeps track of the matching $M$ in a rather cryptic manner; it is not even clear that $m(M)$ determines $M$, although this will follow from a later result (Proposition \ref{facedetermine}). It is natural to add new variables that code directly for the presence or absence of individual edges. Let $i$ and $j$ be integers such that $i+j \equiv 1 \mod 2$. Draw the plane such that $i$ increases in the east-ward direction and $j$ in the north-ward. Let the edges east, north, west and south of $(i,j)$ be labelled $a(i,j)$, $b(i,j)$, $c(i,j)$ and $d(i,j)$ respectively. Introduce new free variables with the same names as the edges. 

We now define an enhanced version of $f$ by 
$$f(n_0,i_0,j_0)=\sum_M m(M) \prod_{e \in M} e$$
where the sum is still over all matchings of the order $n_0$ Aztec diamond. Each matching is now counted by both the monomial $m(M)$ and the product of the edges occuring in that matching. The new $f$ obeys the recurrence
\begin{eqnarray*}
f(n,i,j) &=& \left( a(i+n-1,j)c(i-n+1,j)f(n-1,i,j+1)f(n-1,i,j-1) +  \right. \\
           & & \left. b(i,j+n-1)d(i,j-n+1)f(n-1,i+1,j)f(n-1,i-1,j) \right) / f(n-2,i,j)\
\end{eqnarray*}

\subsection{Somos Sequences and the Octahedron Recurrence} \label{Somos}

Somos has attempted to rebuild the theory of Theta functions on a combinatorial basis. Theta functions are functions of two variables, traditionally called $z$ and $q$. The simplest Theta function is $\Theta(z,q)=\sum_{n=-\infty}^{\infty} e^{2 i n z} q^{n^2}$. From now on, we will concentrate only on the dependence on $z$ and so write $\Theta(z)$. For a general introduction to $\Theta$ functions, see for example chapter XXI of \cite{WhitWat}.

Somos began with a well known result, that, for any $x_0$ and $z_1$, $\Theta(z)$ would obey 
$$\Theta(z_0 + n z_1) \Theta(z_0 + (n-4) z_1) = r \Theta(z_0 + (n-3) z_1) \Theta(z_0 + (n-1) z_1) + s \Theta(z_0 + (n-2) z_1)^2$$
for $r$ and $s$ certain constants depending in a complicated manner on $q$. Somos proposed using this recurrence to compute $\Theta(z_0+n z_1)$ in terms of $r$, $s$ and $\Theta(z_0)$, $\Theta(z_0+z_1)$, $\Theta(z_0+2 z_1)$ and $\Theta(z_0+3 z_1)$. He discovered that $\Theta(z_0+n z_1)$ was a Laurent polynomial in these $6$ variables. This is not difficult to prove by elementary means, but suggests a deeper explanation may be lurking. Moreover, he discovered that all of the coefficients of this Laurent polynomial were positive but was unable to prove this; this will be proven for the first time as a corollary of the results of this paper.

After Somos's work, many other recurrence with surprising Laurentness properties were discovered; see \cite{Gale} for a summary. The family of most importance for this paper is the Three Term Gale-Robinson Recurrence: Let $k$, $a$ and $b$ be positive integers with $a,b < k$. The Three Term Gale-Robinson Recurrence (abbreviated to ``Gale-Robinson Recurrence'' in the remainder of this paper) is
$$g(n) g(n-k) = r g(n-a) g(n-k+a) + s g(n-b) g(n-k+b).$$
Note that the Somos recurrence is the case $(k,a,b)=(4,1,2)$. Note that if $(k,a,b)=(2,1,1)$ and $r=s=1$, we get the recurrence $g(n) g(n-2) = 2 g(n-1)^2$, the recurrence for the number of perfect matchings of the Aztec diamond.

The Gale-Robinson recurrence, for every value of $(k,a,b)$, can be thought of as a special case of the Octahedron recurrence. This observation was first made by Propp and seems to have first appeared in print in \cite{FomZel}. The reduction is performed as follows: suppose that $g$ obeys the Gale-Robinson recurrence. Define $f(n,i,j)$ by $f(n,i,j)=g(\frac{kn+(2a-k)i+(2b-k)j}{2})$. Then $f$ obeys
$$f(n,i,j) f(n-2,i,j) = r f(n-1,i-1,j) f(n-1,i+1,j) + s f(n-1,i,j-1) f(n-1,i,j+1).$$

Thus, if we choose $a(i,j)$, $b(i,j)$, $c(i,j)$ and $d(i,j)$ to be constants $a$, $b$, $c$ and $d$ with $ac=r$ and $bd=s$, $f$ obeys the octahedron recurrence. In these circumstances, we must study the octahedron recurrence not with the initial conditions $f(0,i,j)$ and $f(-1,i,j)$ but with the initial coniditions $f(n,i,j)$ for $(n,i,j)$ such that $-k < \frac{kn+(2a-k) i + (2b-k) j}{2} \leq 0$. 

It therefore beomes our goal to study the octahedron recurrence with general initial conditions. The main result of this paper will be that, with any inital conditions, $f(n,i,j)$ is a Laurent polynomial in the initial coniditions and in the coefficients $a(i,j)$, $b(i,j)$, $c(i,j)$ and $d(i,j)$. Moreover, all of the coefficients of this Laurent Polynomial are $1$. In particular, they are positive. In the reduction of the Gale-Robinson recurrence to the octahedron recurrence, we imposed that many of the initial conditions of the octahedron recurrence be equal to each other. Thus, many monomials that would be distinct if the octahedron recurrence were run with arbitrary inital conditions become the same in the Gale-Robinson recurrence. However, we can still deduce that all of the coefficients of the Laurent polynomials computed by the Gale-Robinson recurrence are positive.

\section{Introduction and Terminology} \label{Intro}

In this section we will introduce the main objects and results of this paper. The main object of study of this paper is a function $f$ defined on a three dimensional lattice by a certain recurrence. We denote the set of initial conditions from which $f$ is generated by $\II$ and consider varying the shape of $\II$ inside this lattice. It turns out that the values of $f$ are Laurent polynomials in the initial variables where every term has coefficient 1. We are able further more to find families of graphs such that $f$ gives generating functions for the perfect matchings of these graphs. Special cases of this result include being able to choose $\II$ so that $f$ gives us the values of three term Gale-Robinson sequences or such that the families of graphs are Aztec Diamonds, Fortresses and other common families of graphs with simple formulas for the number of their matchings. 

In the next subsection, we introduce the vocabulary necessary to define $f$. In the following subsection, we introduce the vocabulary necessary to define the families of graphs.

\subsection{The Recurrence}

Set
$$\LL = \{ (n,i,j) \in \ZZ^3,\ n=i+j \mod 2 \}$$
and set
\begin{eqnarray*}
\EE &=& \{ (i,j,q) \in  \ZZ^2 \times \{a,b,c,d\},\ i+j=1 \mod 2 \} \\ 
\FF &=& \{ (i,j) \in \ZZ^2 \}
\end{eqnarray*} 
where $q$ denotes one of the four symbols $a$, $b$, $c$ and $d$. $\LL$, $\EE$ and $\FF$ stand for ``lattice'', ``edges'' and ``faces.''

We call $(i_1, j_1)$ and $(i_2, j_2) \in \FF$ ``lattice-adjacent'' if  $|(i_1-i_2)|+|(j_1-j_2)|=1$. (Another notion of adjacency will arise later.)

For $(n_0,i_0,j_0) \in \LL$, set 
$$p_{(n_0,i_0,j_0)}(i,j)=n_0-|i-i_0|-|j-j_0|.$$
(We will often drop the subscript on $p$ when it is clear from context.) Let
\begin{eqnarray*}
C_{(n_0,i_0,j_0)} &=& \{ (n,i,j) \in \LL : n \leq p_{(n_0,i_0,j_0)}(i,j) \} \\
\Copen_{(n_0,i_0,j_0)} &=& \{ (n,i,j) \in \LL : n < p_{(n_0,i_0,j_0)}(i,j) \} \\
\Cbound &=& \{ (n,i,j) \in \LL : n = p_{(n_0,i_0,j_0)}(i,j) \}
\end{eqnarray*}
We call $C$, $\Copen$ and $\Cbound$ the \emph{cone}, \emph{inner cone} and \emph{outer cone} of $(n_0,i_0,j_0)$ respectively. $C$ is a square pyramid with it's vertex at $(n_0,i_0,j_0)$.

Let $h \maps \FF \to \ZZ$ and define
\begin{eqnarray*}
\II &=& \{ (i,j,n) \in \LL,\ n = h(i,j) \} \\
\UU &=& \{ (i,j,n) \in \LL,\ n > h(i,j) \}.
\end{eqnarray*}

Assume that
\begin{enumerate}
\item $|h(i_1,j_1)-h(i_2,j_2)|=1$ if $(i_1,j_1)$ and $(i_2, j_2)$ are lattice-adjacent.
\item $(h(i,j),i,j) \in \LL.$
\item $\lim_{|i|+|j| \to \infty} h(i,j)+|i|+|j|=\infty$.
\end{enumerate}

The last condition is equivalent to ``For any $(n,i,j) \in \LL$, $C_{(n,i,j)} \cap \UU$ is finite.'' We will call such an $h$ a ``height function'' and call a function $h$ that obeys the first two conditions a pseudo-height function. Note: this definition is not related to the height functions in the theory of Aztec diamonds, described in \cite{CEP}.

Let $K$ be the field of formal rational functions in the following infinite families of variables: the family $x(i,j)$, $i$, $j \in \ZZ^2$ and the families $a(i,j)$, $b(i,j)$, $c(i,j)$, $d(i,j)$ where $i$, $j \in \ZZ^2$ and $i+j=0 \mod 2$. (Clearly, $\FF$ indexes the $x$'s and $\EE$ indexes the $a$'s, $b$'s, $c$'s and $d$'s). Let $R \subset K$ be the ring $\ZZ[x(i,j), 1/x(i,j), a(i,j), b(i,j), c(i,j), d(i,j)].$ For reasons to appear later, we call the $x$'s the ``face variables'' and the $a$'s, $b$'s, $c$'s and $d$'s the ``edge variables.''

We define a function $f \maps \II \cup \UU \to K$ recursively by $f(n,i,j)=x(i,j)$ for $(n,i,j) \in \II$ and 
\begin{eqnarray*}
f(n,i,j) &=& \left( a(i+n-1,j)c(i-n+1,j)f(n-1,i,j+1)f(n-1,i,j-1) +  \right. \\
           & & \left. b(i,j+n-1)d(i,j-n+1)f(n-1,i+1,j)f(n-1,i-1,j) \right) / f(n-2,i,j) 
\end{eqnarray*}
for $(n,i,j) \in \UU$. Our third condition on $h$ ensures that this recurrence will terminate.

The main result of this paper is
\begin{MT1}
$f(n,i,j) \in R$ and, when written as a Laurent polynomial, each term appears with coefficient 1. Moreover, the exponent of each face variable is between $-1$ and $3$ (inclusive) and the exponent of each edge variable is $0$ or $1$.
\end{MT1}

To describe our result more precisely, we need to introduce some graph-theoretic terminology. 

\subsection{Graphs With Open Faces} \label{gwofs}

Let $G_0$ be a connected planar graph with a specified embedding in the plane and let $e_1$, $e_2$, \dots $e_n$ be the edges of the outer face in cyclic order.  We define a ``graph with open faces'' to be such a graph $G_0$ with a given partition of the cycle $e_1$, \ldots $e_n$ into edge disjoint paths $e_i$, $e_{i+1}$, \dots $e_j$. 

Denote by $G$ the graph with open faces associated to $G_0$ and the partition above. We call a path $e_i$, \ldots $e_j$ in the given partition an ``open face'' of $G$. We denote the open faces of $G$ by $F_o(G)$. We refer an interior face of $G_0$ as a ``closed face'' of $G$ and denote the set of them by $F_c(G)$. We set $F(G)=F_o(G) \cup F_c(G)$ and call a member of $F(G)$ a face of $G$. Note that the exterior face of $G_0$ is \emph{not} usually a face of $G$. 

The image associated to a \gwof is that open disks have had portions of their boundary glued to the outside of $G_0$ along the paths $e_i$, \dots $e_j$ with some additional boundary left hanging off into space. 

We will at times need to allow $G_0$ to be infinite, in which case $G_0$ might have several outer faces, or none at all. We will not describe the appropriate modifications, as they should be obvious. The infinite graphs which will come up in this paper will have no outer face and the corresponding \gwofs will have no open faces.

We will refer to an edge or vertex of $G_0$ as an edge or vertex (respectively) of $G$. We say that an edge $e$ borders an open face $f$ of $G$ if $e$ lies in the associated path $e_i$, \dots $e_j$. All other incidence terminology (endpoint, adjacency of faces, etc.) should be intuitive by analogy to ordinary planar graphs. We use $E(G)$ and $V(G)$ to denote the edges and vertices of $G$.

By a map $G' \hookrightarrow G$ we mean a triple of injections $V(G') \hookrightarrow V(G)$, $E(G') \hookrightarrow E(G)$ and  $F(G') \hookrightarrow F(G)$ compatible with the adjacency relations. We say $G'$ is a sub-\gwof if such a map exists. Note that we do not require that open faces be taken to open faces.

Let $R(G)$ be the ring $\ZZ[x_f, 1/x_f, y_e]$ where $x_f$ are formal variables indexed by $f \in F(G)$ and $y_e$ are formal variables indexed by $e \in E(G)$.

For the rest of this subsection, assume that $G$ is finite. A matching of $G$ is a collection of edges $M$ such that every vertex lies on exactly one edge of $M$. 

Let $M$ be a matching of $G$. For $e\in E(G)$, set $\delta(e)=1$ if $e \in M$, $0$ otherwise. For $f \in F(G)$, let $a$ be the number of edges of $f$ that lie $M$ and $b$ the number of edges of $f$ not in $M$. If $f \in F_c(G)$, set 
$$\epsilon(f)=\left\lceil \frac{b-a}{2} \right\rceil-1 ;$$ 
if $f \in F_o(G)$ set 
$$\epsilon(f)=\left\lceil \frac{b-a}{2} \right\rceil.$$ 
(Note that $\epsilon(f) \geq -1$.)

For any matching $M$, set 
$$m(M)=\prod_e y_e^{\delta(e)} \prod_f x_f^{\epsilon(f)}.$$
For any \gwof $G$, set 
$$m(G)=\sum_M m(M)$$
where the sum is over all matchings of $G$. We call $m(M)$ the matching monomial of $M$ and $m(G)$ the matching polynomial of $G$. 

It will turn out that the actual graphs to which we will apply this definition are bipartite, so every closed face has an even number of edges and we may omit the $\lceil \  \rceil$ in this case. However, it will be convenient to have a definition that is valid for all graphs.

We can now give our second statement of the Main Theorem.

\begin{MT2}
For any height function $h$ we can find an infinite \gwof $\GG$, a decomposition $E(\GG)=E_w \sqcup E_u$ and injective maps $\alpha \maps E_w \to \EE$, $\alpha \maps F \to \FF$ and a family of sub-\gwofs $G_{n,i,j}$ indexed by $(n,i,j) \in \UU$ such that

\begin{enumerate}
\item If we define $\alpha \maps R(\GG) \to R$ by $\alpha(x_f)=x(\alpha(f))$ for $f \in F(\GG)$,  $\alpha(y_e)=q(i,j)$ where $\alpha(e)=(i,j,q)$ for $e \in E_w$ and  $\alpha(y_e)=1$ for $e \in E_u$, then 
$$f(n,i,j)=\alpha(m(G_{n,i,j})).$$
(We have implicitly used the obvious injection $R(G_{n,i,j}) \into R(\GG)$.)
\item If $(n',i',j') \in C_{n,i,j}$ then $G_{n',i',j'}$ is a sub-\gwof of $G_{n,i,j}$
\end{enumerate}
We refer to the edges in $E_w$ as the ``weighted edges'' and the edges in $E_u$ as the ``unweighted edges''.
\end{MT2}

Clearly, this will imply the previous statement of the main theorem, except for the bounds on the exponents and the claim that the coefficients are 1. The first will follow by showing as well that every face of $G$ has $\leq 8$ edges. The second will follow by showing that the edges of $E_u$ are vertex disjoint, so a matching $M$ is uniquely determined by $M \cap E_w$.

\subsection{Plan of the Paper}

In Section~\ref{Algorithm}, we will describe an algorithm we refer to as ``the method of crosses and wrenches'' for finding the graph $\GG$ and the subgraphs $G_{(n,i,j)}$. We will postpone proving the correctness of the algorithm to present applications of the theory in Section~\ref{Examples}. In particular, we will carry out the examples of Somos-4 and Somos-5 and show explicitly the families of graphs associated to them. We will also show that, by choosing certain periodic functions for $h$, we can create fortress graphs and families of graphs studied by Chris Douglass and Matt Blum. The number of matchings of these graphs are powers of 5, 2 and 3 respectively, in each case we will give a rapid proof of this by induction from our Main Theorem. Fortresses and Douglass' graphs are discussed in \cite{Ciucu} and these formulas are also proved there.

In the Section~\ref{Urban}, we will first relate the matchings of the graphs $G_{(n,i,j)}$ to the matchings of certain infinite graphs subject to ``boundary conditions at infinity''. This will remove the elegant property that $G_{(n',i',j')} \hookrightarrow G_{(n,i,j)}$ when $(n',i',j') \in C_{(n,i,j)}$, but it will create objects better suited to an inductive proof.  We will then prove the main theorem by varying $\II$ and holding $(n,i,j)$ fixed. In the process, we will recover variants of the ``Urban Renewal''operations of \cite{Ciucu}.

In Section~\ref{Kuo}, we will give a second proof, inspired by a proof of \cite{Kuo} for the Aztec Diamond case, that works by holding $\II$ fixed and varying $(n,i,j)$. This proof has some additional consequences that the first proof does not. However, due to the extraordinary number of cases that would otherwise be involved, we will not check that all of the exponents work out.

In Section~\ref{Final} we make some final comments.

\subsection{A Note on Simplification of Results}

One might wonder whether the introduction of so many families of variables is necessary to the paper. For most parts of the paper, one may think of any family variables one does not want to deal with as simply being set equal to $1$ without difficulty. However, there is one place where this does not work: our first proof of the main theorem is inductive and the induction only goes through if the theorem is stated with the entire family of face variables present.

\subsection{Acknowledgments}

Most of this research was done while I was working in REACH (Research Experiences in Algebraic Combinatorics at Harvard) and was aided by many helpful conversations with the other members of REACH, and in particular by Daniel Abramson, Gabriel Carroll and Seth Kleinerman's willingness to work through many special cases with me. Above all, I would like to thank Jim Propp for suggesting this problem and showing me much of his unpublished work on it. 

\section{The Method of ``Crosses and Wrenches''} \label{Algorithm}

In this section we describe how to find the graphs $\GG$ and $G_{(n,i,j)}$ discussed above. As a running example, we we'll compute $f(3,1,0)$ with $h(i,j)=|i+j|$. The goal is to predict the formula:
\begin{eqnarray*}
f(3,1,0) &=& a(3,0) c(-1,0) a(2,-1) c(0,-1) x(1,-2) x(1,-1)^{-1} x(1,1) \\
&+& a(3,0) c(-1,0) b(1,0) d(1,-2) x(1,0)^{-1} x(0,-1) x(2,-1) x(1,-1)^{-1} x(1,1) \\
&+& b(1,2) d(1,-2) a(1,0) c(-1,0) x(0,-1) x(0,1) x(0,0)^{-1} x(2,0) x(1,0)^{-1} \\
&+& b(1,2) d(1,-2) b(0,1) d(0,-1) x(-1,0) x(0,0)^{-1} x(2,0).
\end{eqnarray*}

\begin{figure}
\centerline{\scalebox{0.7}{\includegraphics{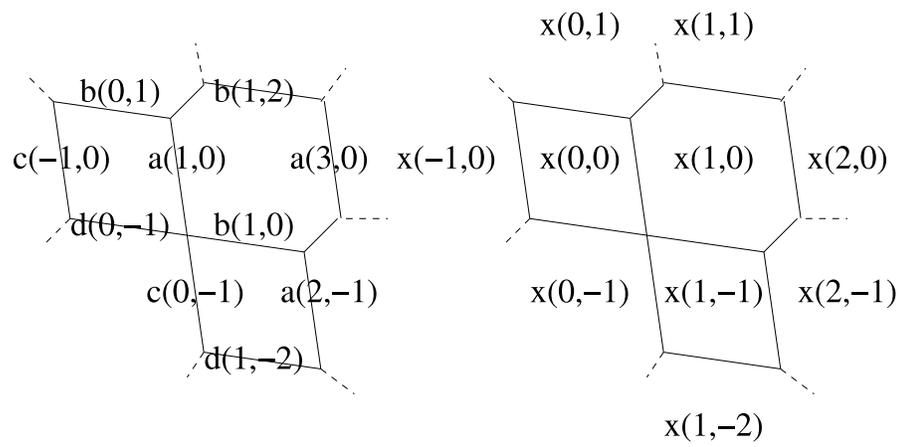}}}
\caption{The Graph of our Running Example} \label{example}
\end{figure}

\begin{figure}
\centerline{\scalebox{0.5}{\includegraphics{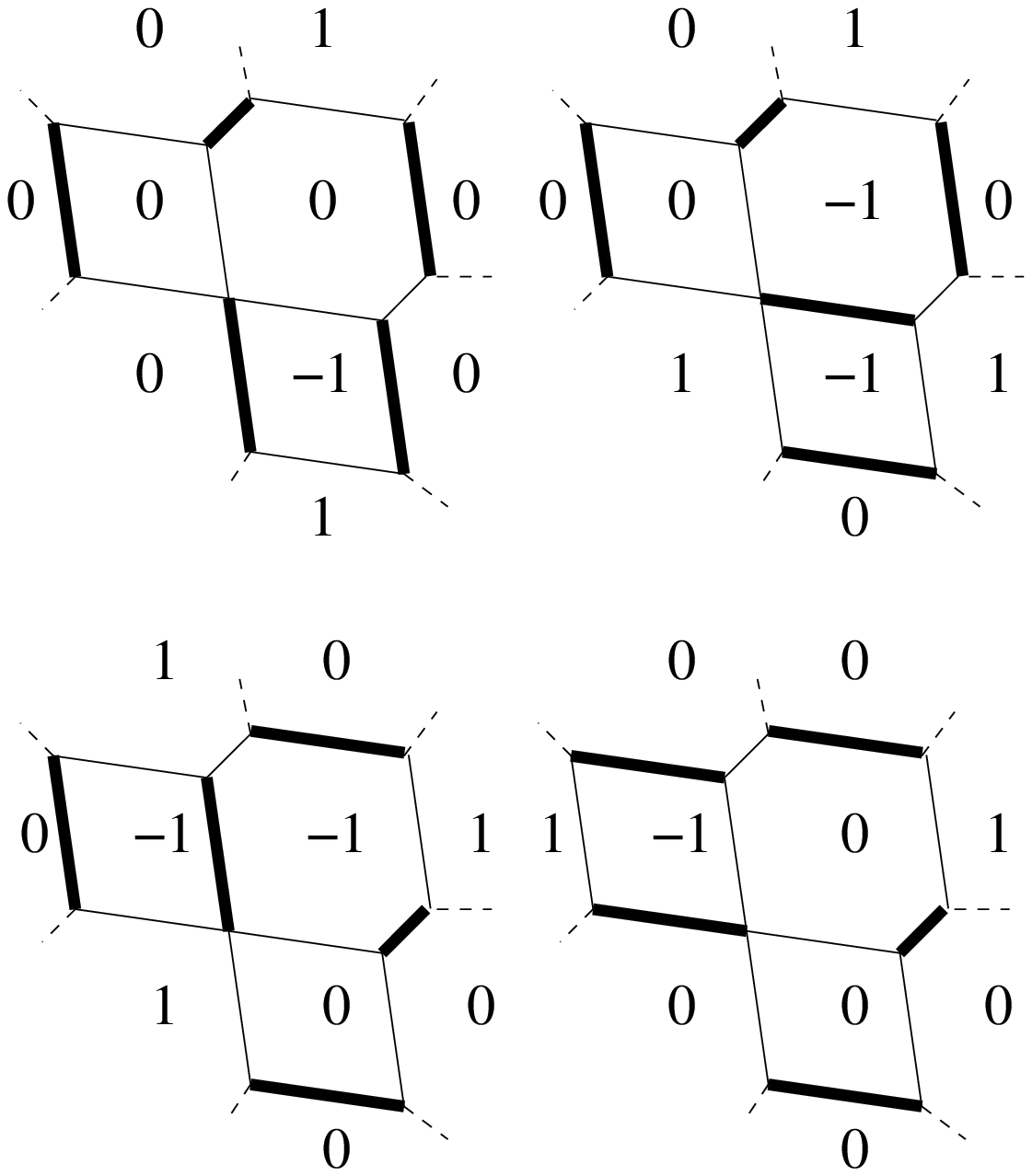}}}
\caption{The Four Matchings of our Running Example} \label{matchings}
\end{figure}

The four terms in the above formula will eventually be shown to correspond to the four matchings of the \gwof shown in Figure \ref{example}. Here the dashed edges seperate the various open faces and the graph is drawn twice so that the face and edge labels will not overlap each other. The four matchings are shown in Figure \ref{matchings}; the numbers in Figure \ref{matchings} are the exponents of the corresponding face variables.

The octahedron recurrence with these initial conditions appears section 9.4 of \cite{FG}. A ``tropicalized'' version of our running example appears in \cite{KTW}.

\subsection{The Infinite Graph $\GG$}

Let $h(i,j)$ be a height function. We describe the graph $\GG$ by describing its dual. The faces of $\GG$ (which are all closed) are indexed by the elements of $\II$ and the map $\alpha$ sends the face $(n,i,j) \in \II$ to $(i,j) \in \FF$. Picture the face $(n,i,j)$ as centered at the point $(i,j) \in \RR^2$. 

If $(i_1, j_1)$ and $(i_2, j_2)$ are lattice-adjacent (i.e. $|(i_1-i_2)-(j_1-j_2)|=1$) then $(n_1, i_1, j_2)$ borders $(n_2, i_2, j_2)$. In addition, if $|i_1-i_2|=|j_1-j_2|=1$ then $(n_1, i_1, j_1)$ borders $(n_2, i_2, j_2)$ iff $h(i_1, j_2) \neq h(i_2, j_1)$. No other pairs of faces border and faces which border border only along a single edge.

We refer to this method of finding $\GG$ as the ``method of crosses and wrenches'' because it can be describe geometrically by the following procedure: any quadruple of values
$$\begin{pmatrix}
h(i,j) & h(i+1, j) \\
h(i,j+1) & h(i+1,j+1)
\end{pmatrix}$$
must be of one of the following six types:

\begin{alignat*}{10}
{\begin{pmatrix}
h & h+1 \\
h+1 & h
\end{pmatrix}} \quad &
{\begin{pmatrix}
h+1 & h \\
h & h+1
\end{pmatrix}} \quad &
{\begin{pmatrix}
h & h+1 \\
h+1 & h+2
\end{pmatrix}} \\
{\begin{pmatrix}
h+2 & h+1 \\
h+1 & h
\end{pmatrix}} \quad &
{\begin{pmatrix}
h+1 & h \\
h+2 & h+1
\end{pmatrix}} \quad &
{\begin{pmatrix}
h+1 & h+2 \\
h & h+1
\end{pmatrix}}
\end{alignat*}

In the center of these squares, we draw a \cross, \backwrench or \wrench in the first and second, third and fourth, and fifth and sixth cases respectively. The diagram below displays the possible cases.

\begin{alignat*}{10}
{\begin{pmatrix}
h & & h+1 \\
& \cross & \\
h+1 & & h 
\end{pmatrix}} \quad &
{\begin{pmatrix}
h+1 & & h \\
& \cross & \\
h & & h+1 
\end{pmatrix}} \quad & 
{\begin{pmatrix}
h & & h+1 \\
& \backwrench & \\
h+1 & & h+2
\end{pmatrix}} \\
{\begin{pmatrix}
h+2 & & h+1 \\
& \backwrench & \\
h+1 & & h
\end{pmatrix}} \quad & 
{\begin{pmatrix}
h+1 & & h \\
& \wrench & \\
h+2 & & h+1
\end{pmatrix}} \quad &
{\begin{pmatrix}
h+1 & & h+2 \\
& \wrench & \\
h & & h+1
\end{pmatrix}}
\end{alignat*}

We then connect the four points protruding from these symbols by horizontal and vertical edges. (We often have to bend the edges slightly to make this work. Kinks introduced in this way are not meant to be vertices, all the vertices come from the center of a $\cross$ or from the two vertices at the center of a $\wrench$ or $\backwrench$.) We refer to the $\cross$'s as ``crosses'' and the $\wrench$'s and $\backwrench$'s as ``wrenches''.

\begin{figure}
\centerline{\includegraphics{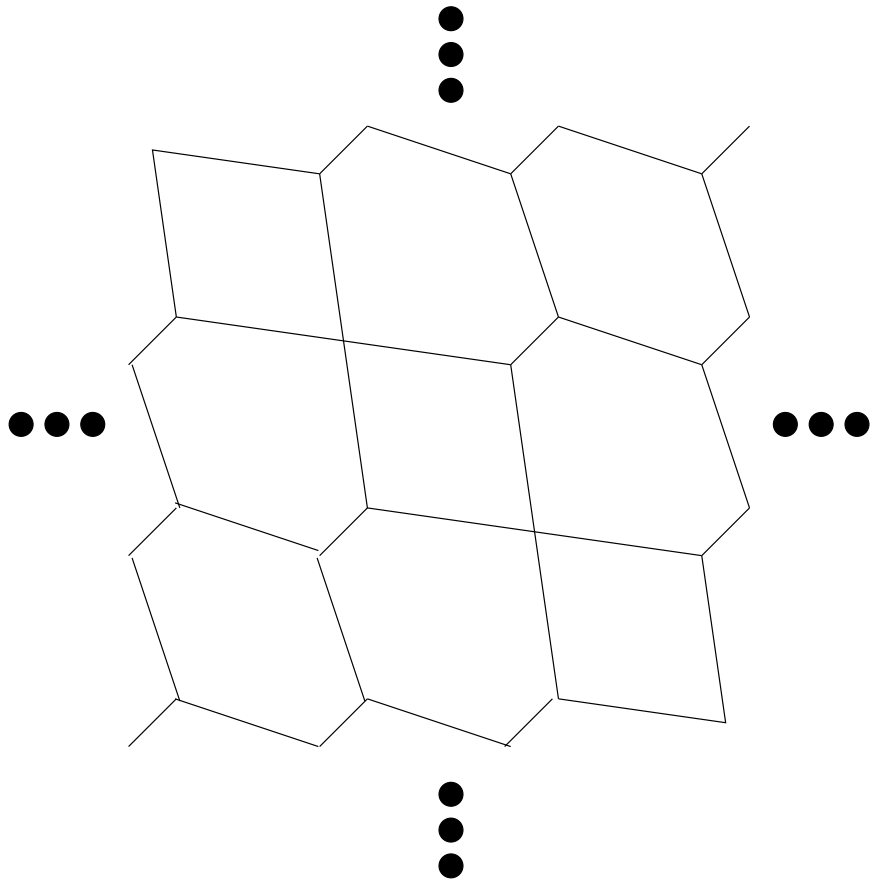}}
\caption{The Infinite Graph of our Running Example} \label{infiniteexample}\end{figure}

In our running example, if $i+j=0$ then $\begin{smallmatrix} h(i,j) & h(i+1, j) \\ h(i,j+1) & h(i+1,j+1) \end{smallmatrix} = \begin{smallmatrix} 0 & 1 \\ 1 & 0 \end{smallmatrix}$ and we place a $\cross$. Otherwise, we get a $\wrench$. the resulting infinite graph, shown in Figure \ref{infiniteexample}, consists of a diagonal row of quadrilaterals seperating a plane of hexagons.

\subsection{Labeling the Edges} \label{LabelEdges}

We will now describe the decomposition $E=E_w \sqcup E_u$ and the map $\alpha \maps E_w \to \EE$.

The set $E_w$ will consist of the horizontal and vertical edges, i.e., those separating lattice-adjacent faces. The set $E_u$ will consist of the diagonal edges, i.e., those which come from the center of a wrench.

Consider any edge $e$ of $E_w$. Such an edge lies between two faces $(n_1, i_1, j_1)$ and $(n_2, i_2, j_2)$ with $|i_1-i_2|+|j_1-j_2|=1$. Without loss of generality, $(i_1-i_2)+(j_1-j_2)=1$. (Otherwise, switch the names of $(i_1, j_1)$ and $(i_2,j_2)$.) There are four cases: 
\begin{enumerate}
\item If $i_1>i_2$ and $n_1>n_2$ then $\alpha(e)=(i_1+n_2,j_1,a)=(i_2+n_1,j_2,a)$.
\item If $i_1>i_2$ and $n_1<n_2$ then $\alpha(e)=(i_1-n_2,j_1,c)=(i_2-n_1,j_2,c)$.
\item If $j_1>j_2$ and $n_1>n_2$ then $\alpha(e)=(i_1,j_1+n_2,b)=(i_2,j_2+n_1,b)$.
\item If $j_1>j_2$ and $n_1<n_2$ then $\alpha(e)=(i_1,j_1-n_2,a)=(i_2,j_2-n_1,d)$.
\end{enumerate}

The reader may wish to refer to Figure \ref{example}.

\subsection{Finding the Subgraphs $G_{(n,i,j)}$}

Finally, we must describe the sub-\gwof $G_{(n,i,j)}$ of $\GG$ that corresponds to a particular $(n,i,j)$. We will abbreviate $G_{(n,i,j)}$ by $G$ in this paragraph. The closed faces of $G$ will be $\Copen_{(n,i,j)} \cap \II$. The edges of $G$ will be the edges of $\GG$ adjacent to some face in $\II \cap G$. The open faces of $G$ will be those members of $\II$ some but not all of whose edges are edges of $G$. Note that every open face of $G$ lies in $\Cbound_{(n,i,j)} \cap \II$ but the converse does not necessarily hold. Note also that there are never any edges that separate two open faces; even if those two open faces are lattice-adjacent. 

Another way to describe the faces of $G$ which is more intuitive but harder to write down is that the faces of $G$ correspond to those $f(n',i',j')$ which are used in compuing $f(n,i,j)$ and the closed faces correspond to those $f(n',i',j')$ which are divided by in the course of this computation.

In our running example, the closed faces are $(0,0,0)$, $(1,1,0)$ and $(0,1,-1)$. The open faces are $(1,0,1)$, $(2,1,1)$, $(1,-1,0)$, $(2,2,0)$, $(1,0,-1)$, $(1,2,-1)$ and $x(1,1,-2)$. 

\subsection{The Main Theorem}

We can now give the final statement of our Main Theorem.

\begin{MT}
For any height function $h$, define $\GG$, $\alpha$ and $G_{(n,i,j)}$ by the algorithm of the previous sections. Then the second statement of the main theorem holds with regard to these choices.
\end{MT}

\subsection{Some Basic Facts} \label{Basic}

It is easy to check that every closed face of $\GG$ has 4, 6 or 8 sides. (Simply check all possible values for the face's eight neighbors.) Hence all of the $G_{(n,i,j)}$ are bipartite. Moreover, every face has $\leq 8$ edges, as promised above. 

It is also easy to fulfill another promise and check that the unweighted edges are vertex disjoint: they lie in the middle of wrenches and do not border each other.

It is clear that $\alpha \maps F(\GG) \to \FF$ is injective (and in fact, bijective). We now show:

\begin{Proposition}
The map $\alpha \maps E_w \to \EE$ is injective.
\end{Proposition}

\begin{proof}
 Consider an edge $e=(i_0,j_0,a) \in \EE$, the cases of $b$, $c$ and $d$ are extremely similar. Any edge with label $e$ must be between two faces of the form $(n,i_0-n-1,j_0)$ and $(n+1,i_0-n,j_0)$. We must show there is at most one value of $n$ for which $(n, i_0-n-1, j_0)$ and $(n+1,i_0-n,j_0)$ both lie in $\II$.

Suppose, for contradiction, there are two such values: $n$ and $n'$. Without loss of generality, suppose that $n<n'$. Then
$$h(i_0-n',j_0)-h(i_0-n-1,j_0)=(n'+1)-n > n'-n-1=(i_0-n')-(i_0-n-1).$$

But $h(i,j)$ must change by $\pm 1$ when $i$ or $j$ changes by 1, so $h$ can not change by this much between $(i_0-n',j_0)$ and $(i_0-n-1,j_0)$, a contradiction.
\end{proof}

We now describe where the faces of $G_{(n_0,i_0,j_0)}$ lie in the plane. We make use of the following abbreviations and notations: we shorten $G_{(n_0,i_0,j_0)}$ to $G$ and $\Copen_{(n_0,i_0,j_0)}$ to $\Copen$. Recall the notation $p(i,j)=n_0-|i-i_0|-|j-j_0|$. So $F_c(G) = \{(i,j)\ :\ h(i,j)<p(i,j) \}$.

\begin{Proposition}
Let $(i,j) \in \FF$ be a closed face of $G$ and let $(i',j') \in \FF$ be such that $i'$ is between $i$ and $i_0$ and $j'$ is between $j$ and $j_0$. Then $(i', j')$ is also a closed face of $G$.
\end{Proposition}

\begin{proof}
Without loss of generality, we may assume that $i_0 \leq i' \leq i$ and $j_0 \leq j' \leq j$. We have $p(i',j')=p(i,j)+(i-i')+(j-j')$. On the other hand, we have $h(i',j') \leq h(i,j) + (i-i')+(j-j')$. As $h(i,j) < p(i,j)$, we have $h(i',j') < p(i',j')$.
\end{proof}

\begin{figure}
\centerline{\scalebox{0.7}{\includegraphics{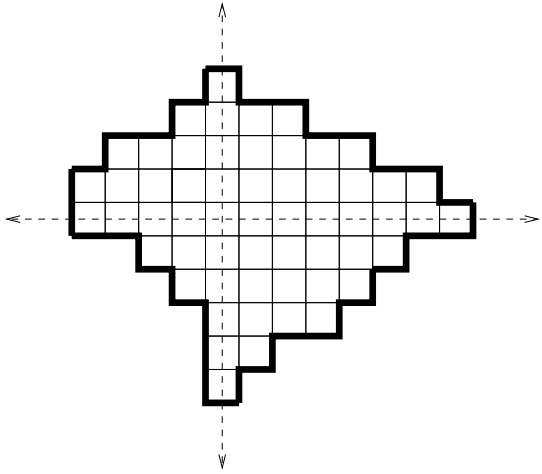}}}
\caption{Positions of Faces of a $G_{(n,i,j)}$} \label{FaceLayout}
\end{figure}

So the faces of $G$ form four ``staircases'', as in figure~\ref{FaceLayout}.

As a corollary, we deduce
\begin{Proposition}
$G$ is connected.
\end{Proposition}

\begin{proof}
Clearly, if $v$ and $w$ lie on the same closed face of $G$, then $v$ and $w$ lie in the same connected component of $G$. But, by the previous proposition, we can travel from any closed face of $G$ to $(i_0, j_0)$ along a sequence of lattice-adjacent faces and lattice adjacent faces have a vertex in common. Thus, since every vertex of $G$ lies on a closed face, every vertex must lie in the same connected component as the vertices of $(i_0,j_0)$. So $G$ is connected.
\end{proof}

It is also clear from the picture above that $G$ is bounded by a single closed loop. Call this loop $S$. Divide $S$, as in figure~\ref{FaceLayout}, into four arcs by horizontal and vertical lines through $(i_0, j_0)$. (In the diagram, $S$ is in bold and the lines through $(i_0,j_0)$ are dashed.)

\begin{Proposition} \label{Boundary}
These lines divide $S$ into four paths. Each of these paths contains an odd number of vertices. On each of these paths, the vertices alternate between vertices all of whose neighbors in $\GG$ are inside or on $S$ (and thus in $G$) and vertices all of whose neighbors in $\GG$ are outside or on $S$. The ends of each path are of the latter type.
\end{Proposition}

\begin{proof}
We consider the section of $S$ on which $i\geq i_0$ and $j \geq j_0$, the other four sections are similar. For the faces $(i,j)$ of $\GG$ within $S$ (that is the closed faces of $G$) we have $h(i,j)<n_0-(i-i_0)-(j-j_0)$ where as for the faces which border $S$ on the outside (that is, the open faces of $G$) we have $h(i,j)=n_0-(i-i_0)-(j-j_0)$. 

Our proof is by induction on the number of faces $(i,j)$ of $G$ for $i \geq i_0$ and $j \geq j_0$. In the base case, where there is only one such face, the path in question has one vertex and the result is trivial. When there is more then one such face, we can always find a face $(i,j)$ such that $(i+1,j)$ and $(i,j+1)$ are outside $S$. For simplicity, we describe the situation where $i>i_0$ and $j>j_0$ and leave the boundary cases to the reader.

Set $h=n_0-(i-i_0)-(j-j_0)$. We must have $h(i+1,j)=h(i,j+1)=h-1$ as they are open faces of $G$. We must have $h(i,j)<h$; as $|h(i,j) - h(i+1,j)|=1$ we must have $h(i,j)=h-2$. If we then change $h(i,j)$ to $h$, we remove $(i,j)$ from within $S$. 

There are four cases, based on whether $h(i-1,j+1)$ and $h(i+1,j-1)=h$ or $h-2$; figure~\ref{BoundaryChange} shows what happens in each case. The edges of $S$ are drawn in bold, the faces of $G$ are surrounded by thin boxes and the bipartite coloring of $G$ is shown in black and white.

\begin{figure}
\centerline{\scalebox{0.5}{\includegraphics{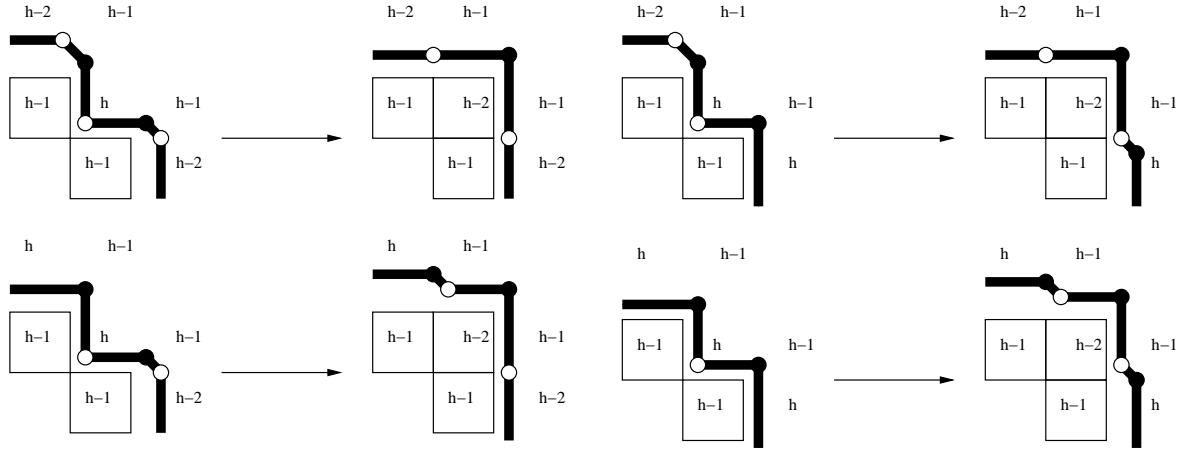}}}
\caption{The Effect on the Boundary of $G_{(n,i,j)}$ of Changing $h(i,j)$}\label{BoundaryChange}
\end{figure}

\end{proof}

We will use this proposition in the proof of proposition \ref{GeomN}. This proposition is also useful in testing whether it is likely that a graph can be realized through the method of crosses and wrenches. For example, the hexagonal graphs of section 6 of \cite{Kuo} can not be so realized, because on their boundary there are six different boundary edges connecting vertices of degree two and a crosses and wrenches graph can only have four such edges.

{Examples} \label{Examples}

\subsection{Aztec Diamonds}

We begin with an example we discussed in the introduction: the Aztec Diamond graphs. In this case, $h(i,j)=0$ or $-1$ when $i+j \mod 2$ is $0$ or $1$ respectively. In this case, every square of four values is 
$$
\begin{matrix} 0 & -1 \\ -1 & 0 \end{matrix} \quad \mathrm{or} \quad 
\begin{matrix} -1 & 0 \\ 0 & -1 \end{matrix}
$$
so in every case we get a cross. 
$$\begin{matrix}
\cross & \cross & \cross & \cross & \cdots \\
\cross & \cross & \cross & \cross & \cdots \\
\cross & \cross & \cross & \cross & \cdots \\
\cross & \cross & \cross & \cross & \cdots \\
\vdots & \vdots & \vdots & \vdots & \ddots
\end{matrix}
$$
The infinite graph $\GG$ is just the regular square grid. The edge labeling is exactly as described in section \ref{AztecIntro}. The graphs $G_{(n,i,j)}$ are the standard Aztec diamond graphs.

\subsection{Fortresses and Douglass' and Blum's Graphs}

In this section we show that Crosses and Wrenches graphs are capable of reproducing several previously studied families of graphs; the number of whose matchings are perfect powers or near perfect powers. Inside each face $(i,j)$, we put the number $h(i,j)$. Most of these families of graphs appear to have first been described in print in \cite{Ciucu}, which we will use as our reference.

Fortress graphs are discussed in \cite{Ciucu}, section four. To obtain a fortress graph, we take 
\begin{alignat*}{10}
h(i,j) &=& &0& \quad i+j &=& 0 \mod 2,\ \phantom{j} & & \\
       &=& &1& \quad i &=& 0 \mod 2,\ j &=& 1 \mod 2 \\
       &=& &-1& \quad i &=& 1 \mod 2,\ j &=& 0 \mod 2.
\end{alignat*}

Every square of values is of the form 
$$
\begin{matrix} 0 & 1 \\ -1 & 0 \end{matrix},\quad 
\begin{matrix} 0 & -1 \\ 1 & 0 \end{matrix},\quad 
\begin{matrix} 1 & 0 \\ 0 & -1 \end{matrix} \quad \mathrm{or} \quad 
\begin{matrix} -1 & 0 \\ 0 & 1 \end{matrix}
$$
so in every case we get a wrench, in the repeating pattern
$$\begin{matrix}
\wrench & \backwrench & \wrench & \backwrench & \cdots \\
\backwrench & \wrench & \backwrench & \wrench & \cdots \\
\wrench & \backwrench & \wrench & \backwrench & \cdots \\
\backwrench & \wrench & \backwrench & \wrench & \cdots \\
\vdots & \vdots & \vdots & \vdots & \ddots
\end{matrix}
$$

\begin{figure}
\centerline{\scalebox{0.5}{\includegraphics{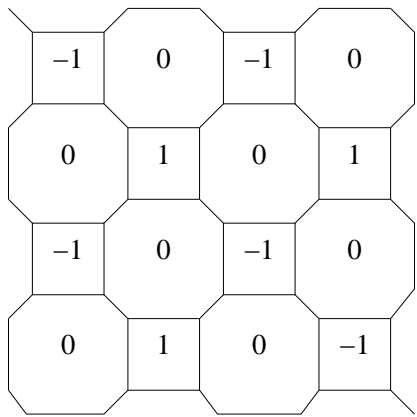}}}
\caption{A Portion of $\GG$ for Fortresses}\label{Fortress}
\end{figure}

Joining up the wrenches, we get the infinite graph in figure~\ref{Fortress}.

Once again, if we set all of the variables equal to 1, $f(n,i,j)$ will count the matchings of fortresses. We must have $f(h(i,j),i,j)=1$ and $f$ must obey the defining recurrence. It is easy to check that both conditions are obeyed by
\begin{eqnarray*}
f(2n,i,j) &=& 5^{n^2} \\
f(2n+1,i,j) &=& 5^{n^2+n} \qquad i=n \mod 2 \\
f(2n+1,i,j) &=& 2 \cdot 5^{n^2+n} \qquad j=n \mod 2,
\end{eqnarray*}
a result of \cite{Ciucu}.

\begin{figure}
\centerline{\scalebox{0.5}{\includegraphics{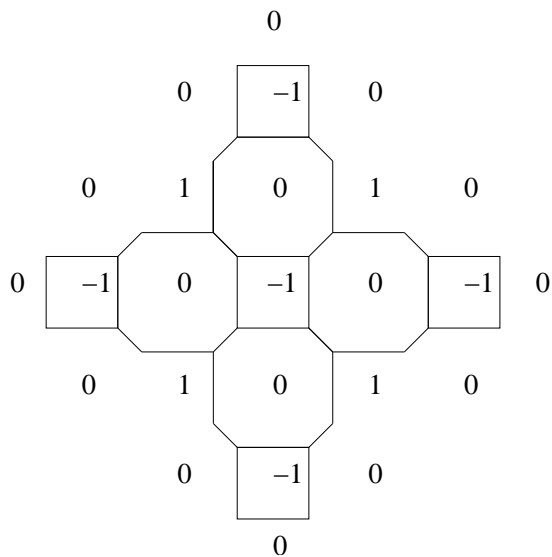}}}
\caption{An Example of a Fortress}\label{FortExamp}
\end{figure}

For example, figure~\ref{FortExamp} shows a fortress with 25 matchings, it is of the form $G_{3,i,j}$ where the center square is $(i,j)$ and $i=0 \mod 2$, $j=1 \mod 2$. This fortress has 9 closed and 16 open faces.

Similarly, to get another family of graphs first studied by Chris Douglas (see \cite{Ciucu}, section six), we take 
\begin{alignat*}{10}
h(i,j) &=& 0 \quad i+j &=& 0 \mod 2, \\
       &=& 1 \quad i+j &=& 1 \mod 4 \\ 
       &=& -1 \quad i+j &=& -1 \mod 4.
\end{alignat*}

We now sometimes will have
$$
\begin{matrix} 1 & 0 \\ 0 & 1 \end{matrix} \quad \mathrm{or} \quad 
\begin{matrix} -1 & 0 \\ 0 & -1 \end{matrix}
$$
which will produce a cross and sometimes will have
$$
\begin{matrix} -1 & 0 \\ 0 & 1 \end{matrix} \quad \mathrm{or} \quad 
\begin{matrix} 1 & 0 \\ 0 & -1 \end{matrix}
$$
which will produce a $\backwrench$.

\begin{figure}
\centerline{\scalebox{0.5}{\includegraphics{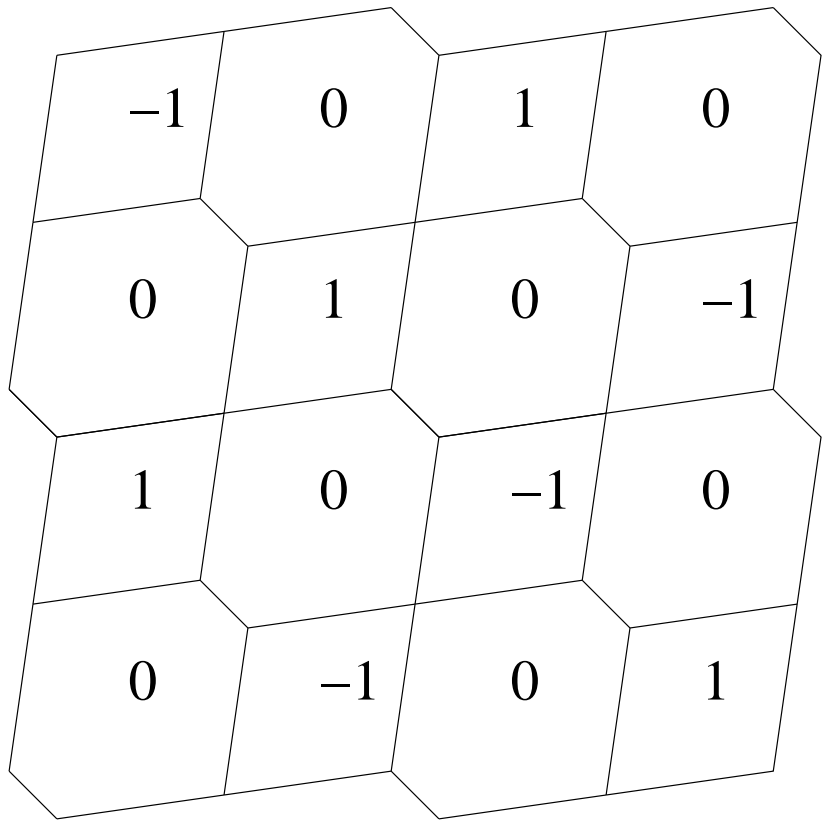}}}
\caption{A Portion of the $\GG$ for Douglass' Graphs}\label{Douglass}
\end{figure}

Overall, we have the repeating pattern
$$\begin{matrix}
\cross & \backwrench & \cross & \backwrench & \cdots \\
\backwrench & \cross & \backwrench & \cross & \cdots \\
\cross & \backwrench & \cross & \backwrench & \cdots \\
\backwrench & \cross & \backwrench & \cross & \cdots \\
\vdots & \vdots & \vdots & \vdots & \ddots
\end{matrix}
$$
which produces the infinite graph in figure~\ref{Douglass}

Again, a quick induction verifies that 
\begin{eqnarray*}
f(2n,i,j) &=& 2^{2 n^2} \\
f(2n+1,i,j) &=& 2^{2 n^2+2 n} \qquad i+j=2n+1 \mod 4 \\
f(2n+1,i,j) &=& 2^{2 n^2+2 n+1} \qquad i+j=2n-1 \mod 4.
\end{eqnarray*}

Similarly, to get the graphs considered by Matthew Blum (see \cite{Ciucu}, section 8), take 
\begin{alignat*}{10}
h(i,j) &=& 0 \quad i+j &=& 0 \mod 2 & & \\
       &=& 1 \quad i+j &=& 1 \mod 2,\ j &=& 0,1 \mod 4 \\ 
       &=& -1 \quad i+j &=& 1 \mod 2,\ j &=& 2,3 \mod 4.
\end{alignat*}
(Ciucu's methods require him to embed $\GG$ into the square grid, and he thus takes $\GG$ to be made of octagons and hexagons. For our purposes, it is more natural to apply Lemma~\ref{split} of Section~\ref{BasicLemmas} to obtain a grid of quadrilaterals and hexagons with the same number of matchings.) We get
$$\begin{matrix}
\wrench & \backwrench & \wrench & \backwrench & \cdots \\
\cross  & \cross      & \cross  & \cross      & \cdots \\
\wrench & \backwrench & \wrench & \backwrench & \cdots \\
\cross  & \cross      & \cross  & \cross      & \cdots \\
\vdots & \vdots & \vdots & \vdots & \ddots
\end{matrix}
$$

\begin{figure}
\centerline{\scalebox{0.5}{\includegraphics{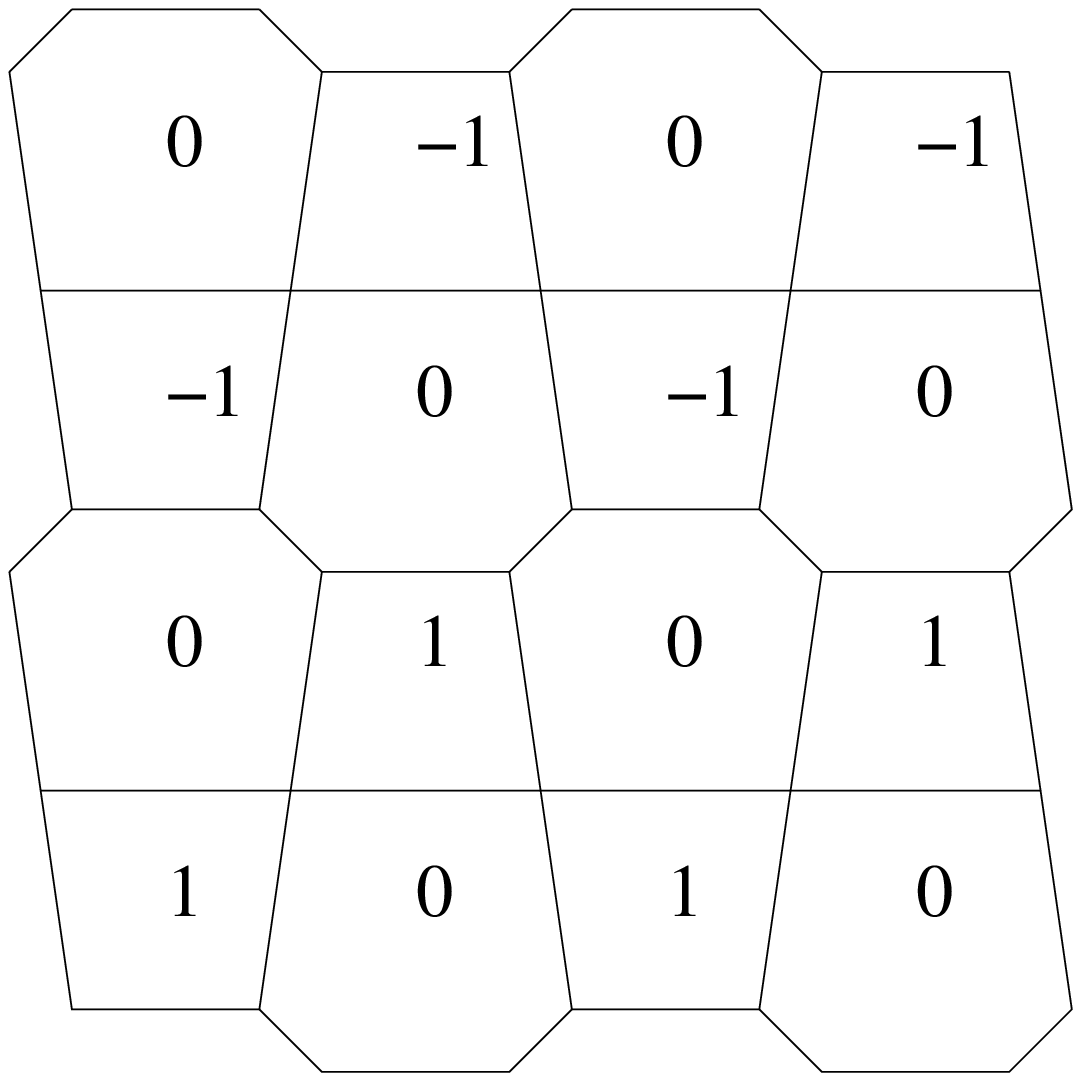}}}
\caption{A Portion of the $\GG$ for Blum's Graphs}\label{Blum}
\end{figure}

This gives us the infinite graph of figure~\ref{Blum}. A quick induction then shows that
\begin{eqnarray*}
f(3n, i, j) &=& \phantom{{} 2 \cdot {}} 3^{3 n^2} \\
f(3n \pm 1, i, j) &=& \phantom{{} 2 \cdot {}} 3^{n^2 \pm n} \qquad \left\lfloor \frac{3n \pm 1}{3} \right\rfloor + \left\lfloor \frac{j}{2} \right\rfloor = 1 \mod 2 \\
f(3n \pm 1, i, j) &=& 2 \cdot 3^{n^2 \pm n} \qquad \left\lfloor \frac{3n \pm 1}{3} \right\rfloor + \left\lfloor \frac{j}{2} \right\rfloor = 0 \mod 2.
\end{eqnarray*}

\subsection{Somos 4 and 5}

As described in the introduction, a combinatorial formula for the octahedron recurrence will give rise to a combinatorial interpretation of the Somos sequences. Specifically, if we set $\II=\{ (n,i,j) \in \LL : -4 < 2n+i \leq 0 \}$ then the number of terms of $f(0,2k,0)$ will obey the Somos-4 recurrence. Similarly, if $\II = \{ (n,i,j) \in \LL : -5 < (5n + i + 3j)/2 \leq 0 \}$ then the number of terms of $f(0,2k,0)$ will obey the Somos-5 recurrence.

\begin{figure}
\centerline{\scalebox{0.5}{\includegraphics{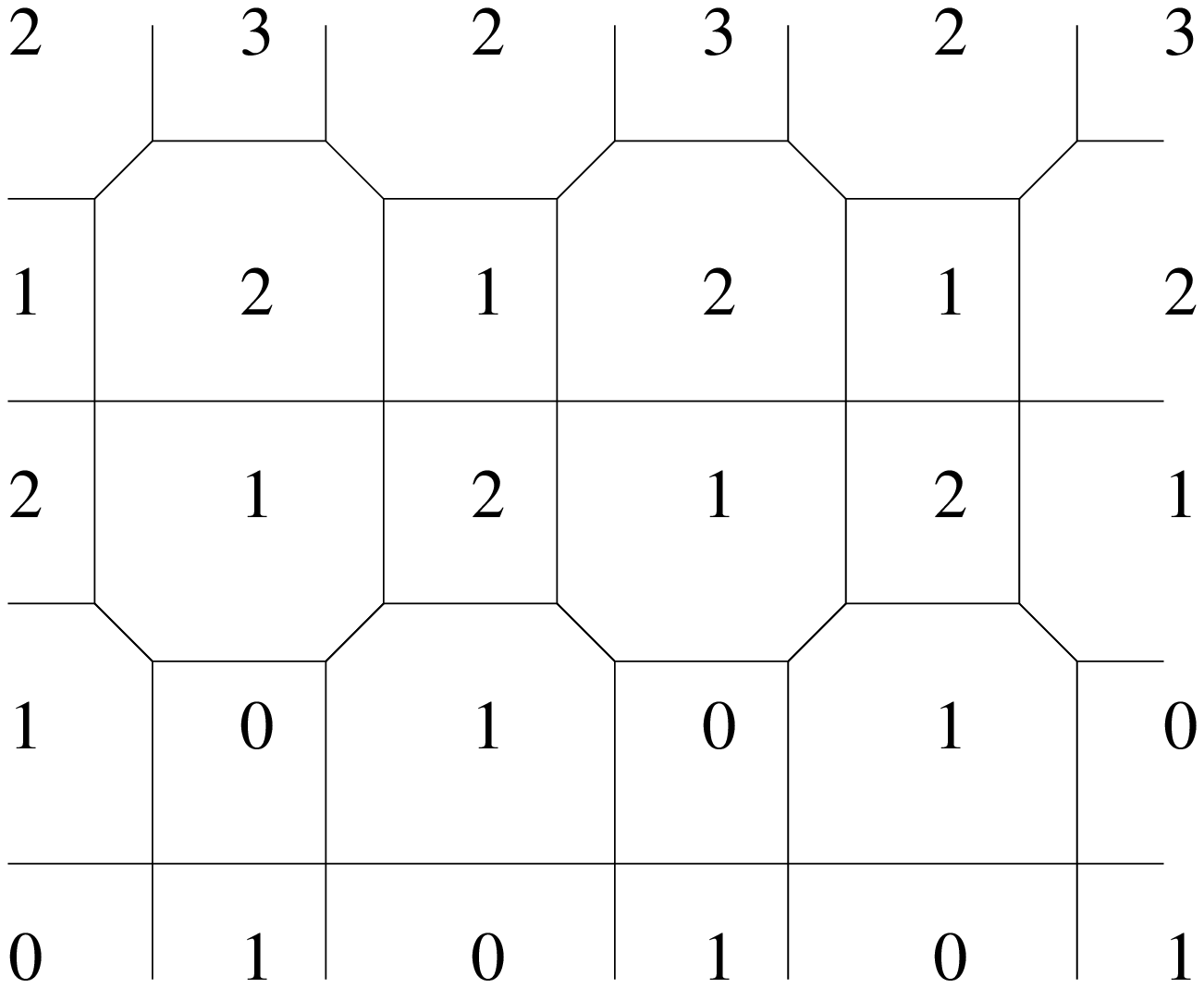}}}
\caption{The $\GG$ for the Somos-4 Sequence}\label{Somos4}
\end{figure}

\begin{figure}
\centerline{\scalebox{0.7}{\includegraphics{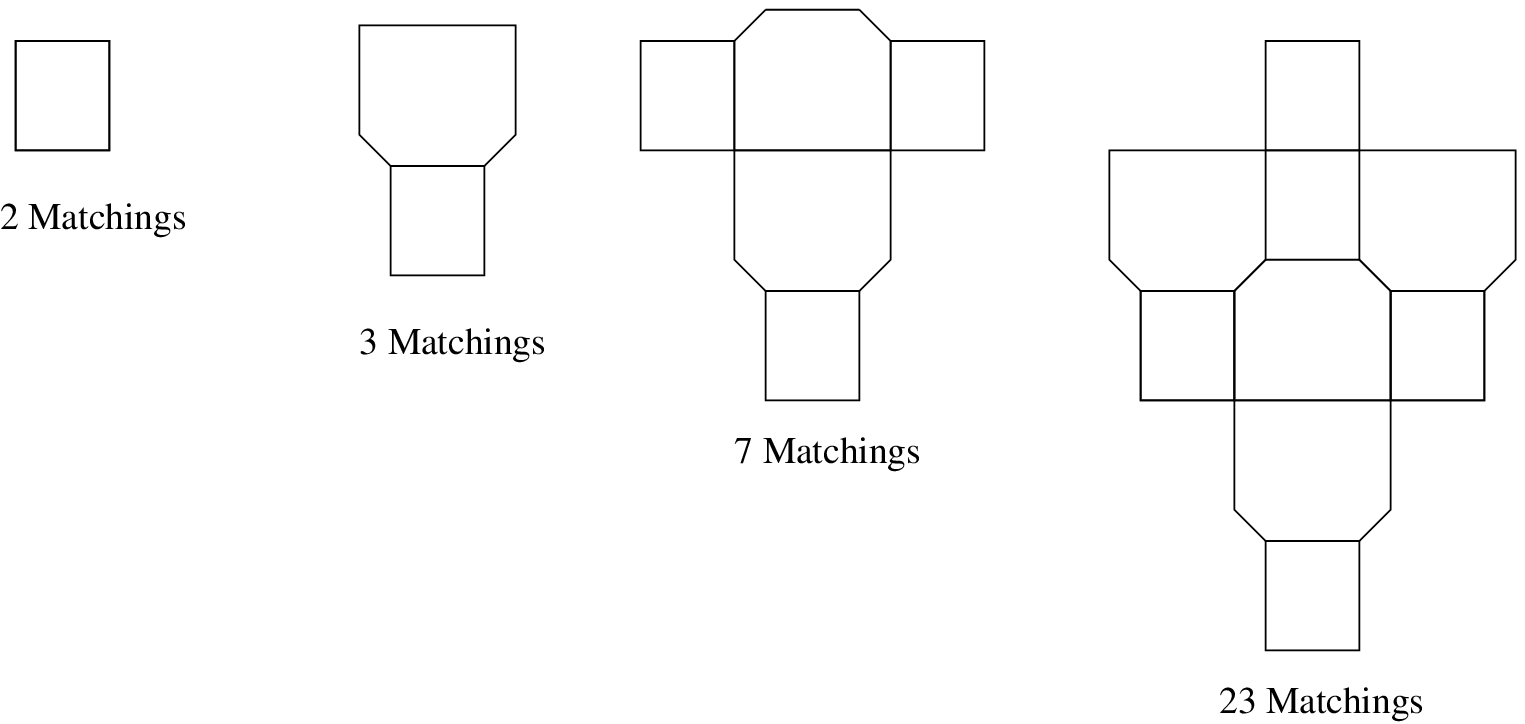}}}
\caption{The First Four Nontrivial Somos-4 Graphs}\label{Somos4Finite}
\end{figure}

\begin{figure}
\centerline{\scalebox{0.5}{\includegraphics{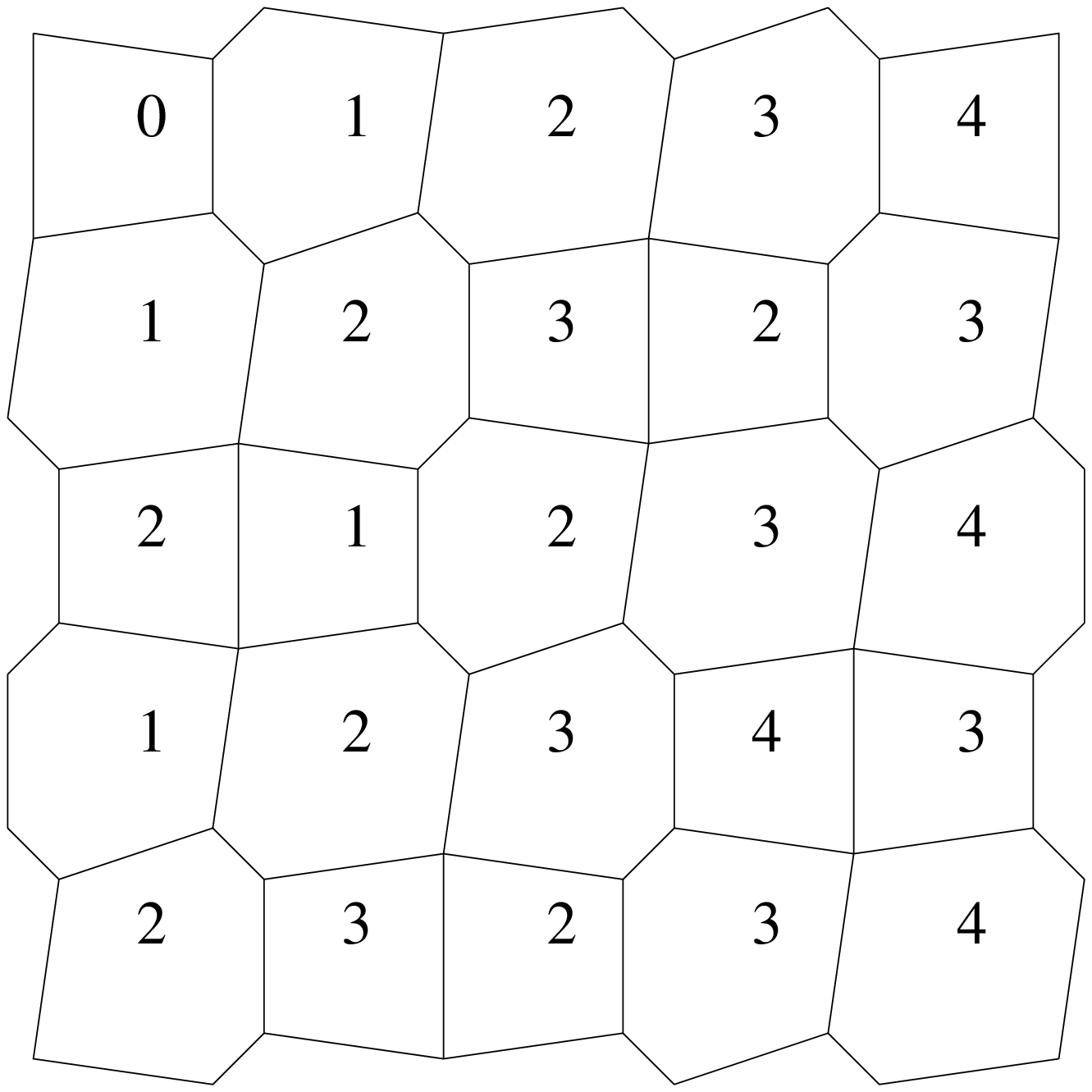}}}
\caption{The $\GG$ for the Somos-5 Sequence}\label{Somos5}
\end{figure}

\begin{figure}
\centerline{\scalebox{0.6}{\includegraphics{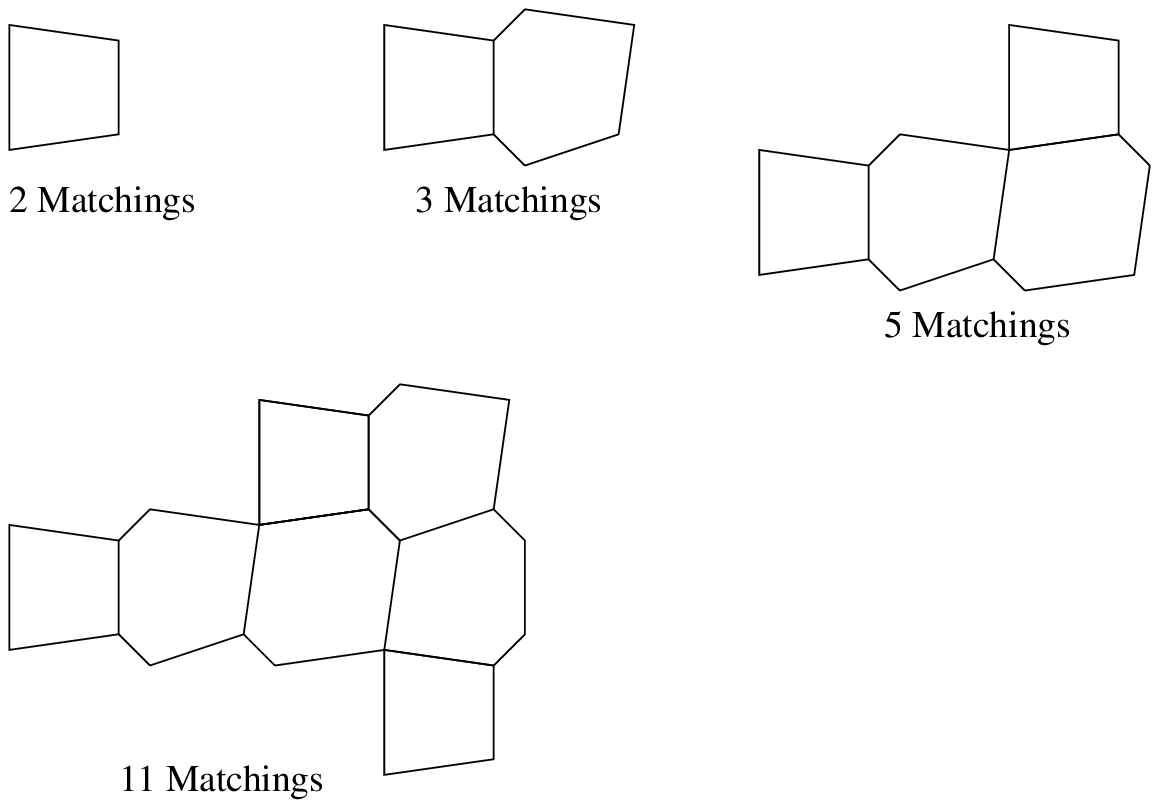}}}
\caption{The First Four Nontrivial Somos-5 Graphs}\label{Somos5Finite}
\end{figure}

Figures~\ref{Somos4} and \ref{Somos4Finite} show the infinite graph and the first several finite graphs corresponding to the Somos-4 recurrence. Figures~\ref{Somos5} and \ref{Somos5Finite} do the same for Somos-5. It may be hard to see the periodicity of the $\GG$ for $(5,1,2)$; the tiles repeat as one travels five cells to the right, or one cell up and two to the right.

\section{Proof I: Urban Renewal} \label{Urban}

In this section, we will prove the correctness of the crosses-and-wrenches algorithm. Our basic strategy is to hold the point $(n_0,i_0,j_0)$ where we are evaluating $f$ fixed, while varying $h$. Our proof is by induction on the number of points in $\UU \cap C_{n_0,i_0,j_0}$. 

\subsection{Infinite Completions} \label{infext}

In this section, we introduce an alternative way of viewing $f(n,i,j)$ as counting the number of matchings of an infinite graph, subject to a ``condition at infinity'' to be described later. This means that \emph{a priori} $f(n,i,j)$ could contain an infinite number of terms, although in fact it will not, and describing the terms of $f(n,i,j)$ requires an \emph{a priori} infinite amount of information. On the other hand, this new method will no longer require the use of open faces and will remove the need to treat faces in the boundary as a special case in our proof. We will use this interpretation in our first proof of the Main Theorem.

Let $h$ be a height function. Fix a triple $(n_0,i_0,j_0) \in \UU$, we will be interested in $f(n_0,i_0,j_0)$. Let $G$ be the \gwof $G_{(n_0,i_0,j_0)}$. Set
$$\tilde{h}(i,j)=\min( h(i,j), p_{(n_0,i_0,j_0)}(i,j)).$$

\begin{figure}
\centerline{\includegraphics{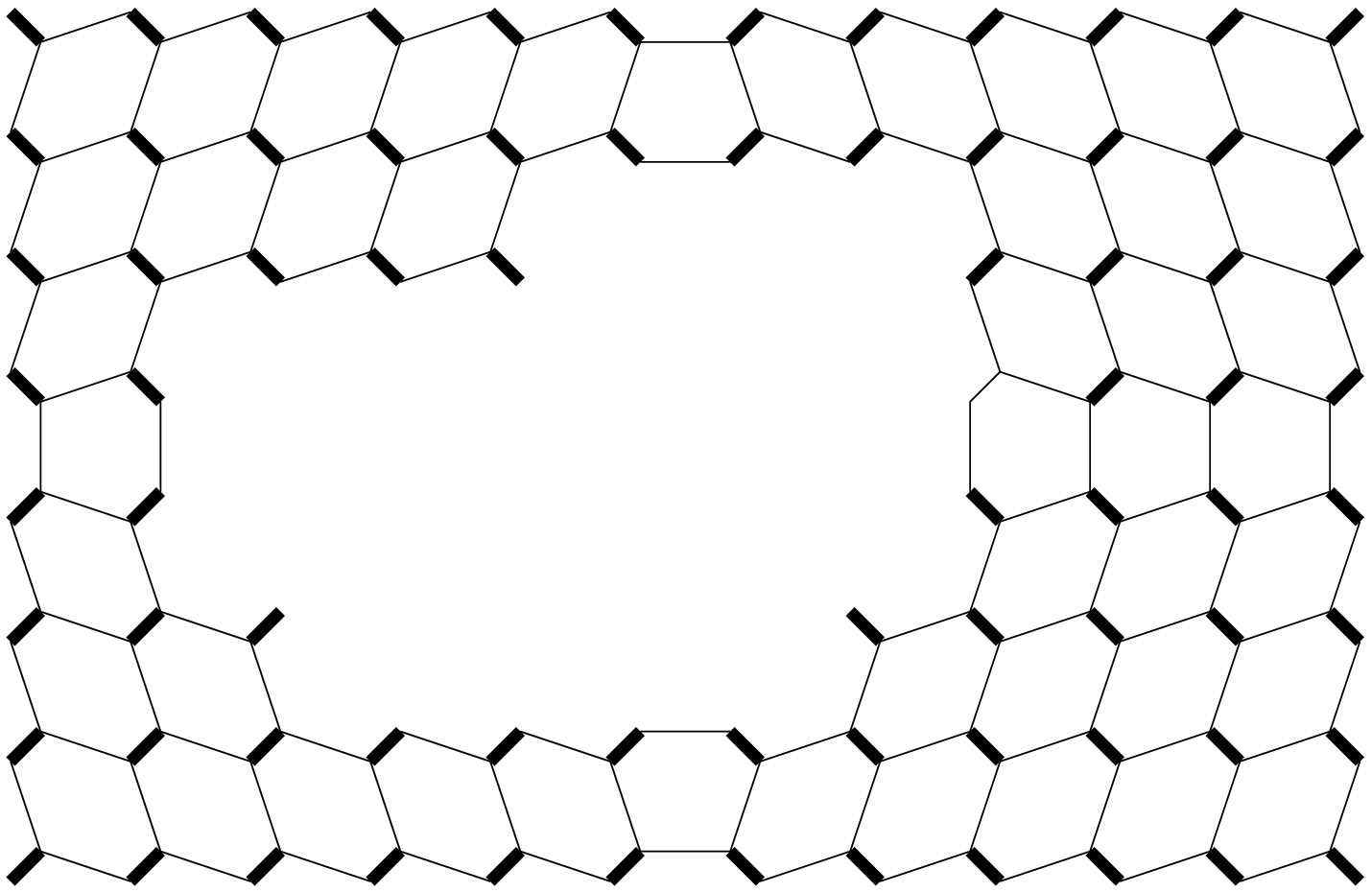}}
\caption{The Standard Matching of $\Gout$}\label{InfComp}
\end{figure}

It is easy to check that $\tilde{h}$ is a pseudo-height function. We can use the method of crosses and wrenches to associate an infinite graph $\tilde{\GG}$ to $\tilde{h}$. Note that if $(n,i,j)$ is a face of $G$ then $h(i,j)=\tilde{h}(i,j)$. Thus, $G$ is a subgraph of $\tilde{\GG}$. Upon removing $G$, what remains looks like figure~\ref{InfComp}. Let $\Gout=\tilde{\GG} \setminus G$.

There is a unique matching (indicated in thick lines) of $\Gout$ given by taking the middle edge of every wrench. Call this matching $\Mout$. Clearly, matchings of $\tilde{\GG}$ that coincide with $\Mout$ on $\Gout$ are in bijection with matchings of $G$. Moreover, one can easily check that the face exponents associated to a given matching of $G$ are the same as for the corresponding matching of $\tilde{\GG}$, including that the exponent of any face of $\tilde{\GG}$ that does not correspond to a face of $G$ is 0.

Thus, we could define $f(n_0,i_0,j_0)$ as the sum over matchings $M$ of $\tilde{\GG}$ such that $M$ coincides with $\Mout$ on $\Gout$ of $m(M)$.

We claim that the following, less obvious result also holds:
\begin{Proposition} Let $M$ be a matching of $\tilde{\GG}$ in which all but finitely many vertices are matched to the same vertex as in $\Mout$. ($\Mout$ is a matching of $\Gout$, so for all but finitely many vertices of $\tilde{\GG}$ this makes sense.) Then $M$ coincides with $\Mout$ everywhere on $\Gout$.
\end{Proposition}

\begin{proof}
Suppose the opposite. All the vertices of $G \setminus \tilde{\GG}$ come from wrenches. Each wrench has one of its vertices closer to $(i_0,j_0)$ than the other; call these the \emph{near} and \emph{far} vertex of the wrench. Now, suppose the theorem is false. Of the finite number of vertices of $\tilde{\GG} \setminus G$ not matched as in $M$, let $v$ be the furthest from $(i_0,j_0)$. We derive a contradiction in the two possible cases.

Case 1: $v$ is a near vertex. Then the far vertex $w$ in the same wrench as $v$ is matched as in $M$. But $w$ is matched with $v$ in $M$, so $w$ is matched with $v$ and $v$ is matched as in $\Mout$.

Case 2: $w$ is a far vertex. But all of $v$'s neighbors except the one it is matched to in $\Mout$ are farther from $(i_0,j_0)$ than $v$ is, so they are matched as in $\Mout$ and can not be matched to $v$. So $v$ is matched as in $\Mout$.
\end{proof}

As result, we can give another statement of the main theorem.

\begin{MTInfExt}
Call a matching $M$ of $\tilde{\GG}$ an acceptable matching if $M$ coincides with $\Mout$ for all but finitely many edges. Then
$$f(n_0,i_0,j_0)=\sum_{M \mathrm{\ acceptable}} m(M).$$
\end{MTInfExt}

\subsection{Some Easy Lemmas}\label{BasicLemmas}

\begin{Lemma} \label{split}
Consider a \gwof $G$. Replace a vertex $G$ with two unweighted edges as in figure~\ref{vertsplit} to form a new $G'$.

\begin{figure}
\centerline{\includegraphics{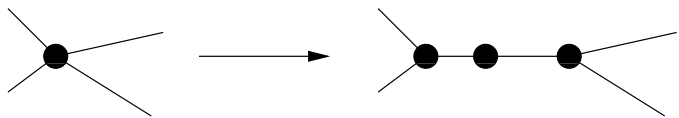}}
\caption{Splitting a Vertex}\label{vertsplit}
\end{figure}

Then $M(G)=M(G')$.
\end{Lemma}

\begin{proof}
Any matching of the old graph can be uniquely transformed to a matching of the new graph by adding one of the two new edges. The additional edge in the matching contributes a factor of 1 in the product of the edge variables. The two faces adjacent to the new edges have one more used edge and one more unused edge than the corresponding faces of the previous matching, so they have the same exponent. The other faces are unaffected.
\end{proof}

%

\begin{figure}
\centerline{\scalebox{0.5}{\includegraphics{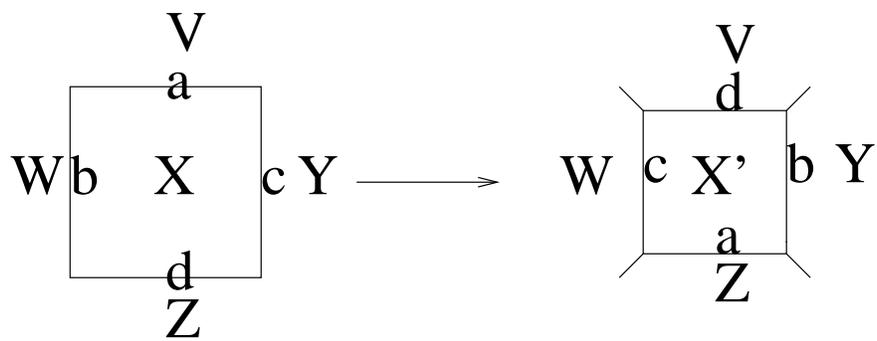}}}
\caption{Applying Urban Renewal}\label{UR}
\end{figure}

\begin{URT}
Suppose a \gwof $G$ contains the sub-\gwof shown in the left of figure~\ref{UR} with the indicated edge weights and face weights. (The face weights are in uppercase, the edge weights in lower case.) Create a new graph $G'$ by replacing this with the \gwof in the right of figure~\ref{UR} where
$$X'=(adWY+bcVZ)/X.$$
Then $m(G')=m(G)$.
\end{URT}

This theorem is related to Lemma 2.5 of \cite{Ciucu}. However, because he doesn't use face weights, Ciucu's theorem involves slightly different replacements and yields a relation of the form $m(G)=(ac+bd) m(G')$.

\begin{proof}
Write $m(G)=m_0+m_1+m_2$ where $m_0$ consists of the terms corresponding to matchings where none of the edges $a$, $b$, $c$ or $d$ is used, $m_1$ consists of the terms using one such edge and $m_2$ consists of the terms using two. Write $m(G')=m'_0+m'_1+m'_2$ similarly. We will show that $m_1=m'_1$, $m_0=m'_2$ and $m_2=m'_0$.

We first give a bijection between the matchings counted by $m_1$ and those counted by $m'_1$ as shown in figure~\ref{biject1}.

\begin{figure}
\centerline{\includegraphics{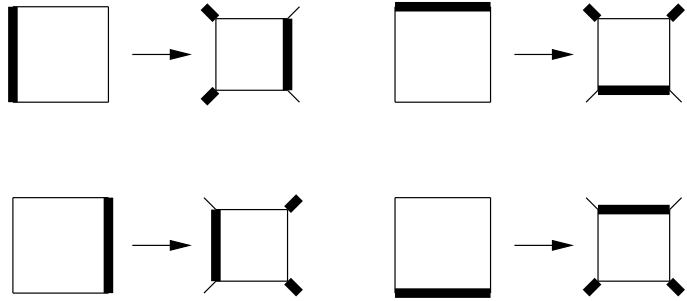}}
\caption{The Bijection in the $m_1$ Case}\label{biject1}
\end{figure}

In each case, it is easy to check that the two matchings paired off have exponent $0$ on $X$ and $X'$ and raise all other variables to the same exponent. 

Next, as shown in figure~\ref{biject0}, we give a bijection that associates to each matching counted by $m_0$ two of the type counted by $m'_2$.

\begin{figure}
\centerline{\includegraphics{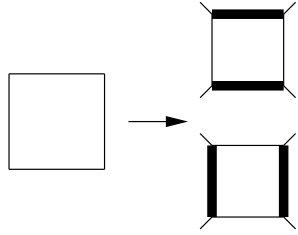}}
\caption{The One-to-Two Map in the $m_0$ Case}\label{biject0}
\end{figure}

The first transformation in figure~\ref{biject0} increases the exponents of $W$ and $Y$ by 1, adds the edges $a$ and $d$ and changes the contribution of the center face from $X$ to ${X'}^{-1}$. Similarly, the second transformation increases the exponents of $V$ and $Z$ by 1, adds the edges $b$ and $c$ and changes the contribution of the center face from $X$ to ${X'}^{-1}$. So we get from $m_0$ to $m'_2$ by multiplying by 
$$\frac{adWY{X'}^{-1}+bcVW{X'}^{-1}}{X}.$$
But, as $X'=(adWY+bcVZ)/X$, this is 1.

Finally, we show $m_2=m'_0$. This is similar to the preceding paragraph. This time, we pair off every two matchings in $m_2$ with one matching in $m'_0$, as shown in figure~\ref{biject2}. 

\begin{figure}
\centerline{\includegraphics{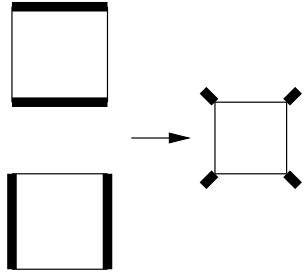}}
\caption{The Two-to-One Map in the $m_2$ Case}\label{biject2}
\end{figure}

In the first pairing of figure~\ref{biject2}, we delete edges $a$ and $d$, reduce the exponents of $W$ and $Y$ by 1 and replace $X^{-1}$ by $X'$. In the second, we delete edges $b$ and $c$, reduce the exponents of $V$ and $Z$ by 1 and replace $X^{-1}$ by $X'$. Thus, we are multiplying by 
$$\frac{X'}{adWYX^{-1}+bcVWX^{-1}}$$
which is 1 as before.

\end{proof}

\subsection{Proof of the Main Theorem}

Let $(n_0, i_0, j_0) \in \UU$ and let $h$ and $\tilde{h}$ be as in subsection \ref{infext}; let $\tilde{\GG}$ be the \gwof associated to $\tilde{h}$. Our proof is by induction on $\sum_{i,j} (p(i,j)-\tilde{h}(i,j))=2 \# (\UU \cap C_{(n_0,i_0,j_0)})$. 

\begin{figure}
\centerline{\includegraphics{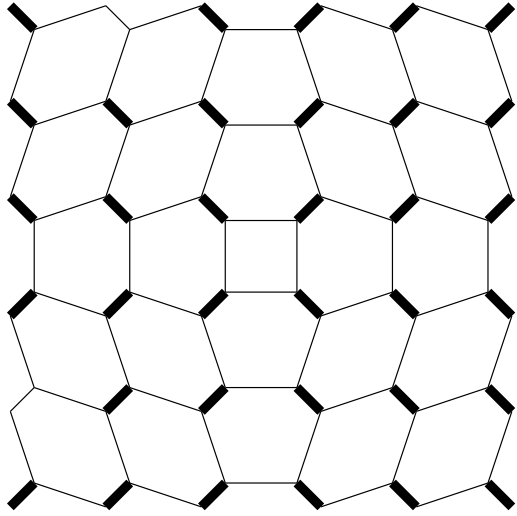}}
\caption{The Base Case, with its Unique Admissible Matching}\label{basecase}
\end{figure}

If $\sum_{i,j} (p(i,j)-\tilde{h}(i,j))=0$ then $\tilde{h}(i,j)=p(i,j)$ so $\tilde{\GG}$ is the graph shown in figure~\ref{basecase} which has only the indicated admissible matching. The matching polynomial of this matching is $x(i_0,j_0)$. We also have in this case that $\tilde{h}(i_0,j_0)=p(i_0,j_0)=n_0$ so we have $f(n_0,i_0,j_0)=x(i_0,j_0)$.

If $\sum_{i,j} p(i,j)-\tilde{h}(i,j) > 0$ then, among the $(i,j)$ for which $h(i,j) < p(i,j)$ there is (at least) one with $\tilde{h}(i,j)$ minimal, let this be $(i,j)$. We claim that $\tilde{h}(i \pm 1,j)=\tilde{h}(i,j \pm 1)=\tilde{h}(i,j)+1$. We deal with the case of $\tilde{h}(i+1,j)$, the other three cases are similar. We have $\tilde{h}(i+1,j)=\tilde{h}(i,j) \pm 1$ as $\tilde{h}$ is a pseudo-height function. If $\tilde{h}(i+1,j)=\tilde{h}(i,j)-1$, then 
$$\tilde{h}(i+1,j) < \tilde{h}(i,j) < p(i,j)-1 \leq p(i+1,j),$$
contradicting the minimality of $(i,j)$.

Thus, we have shown $\tilde{h}(i \pm 1,j)=\tilde{h}(i,j \pm 1)=\tilde{h}(i,j)+1$. So the face associated to $h(i,j)$ is a square. Apply the Urban Renewal Theorem to this square to create a new graph with the same matching polynomial. Then apply Lemma \ref{split} to this graph wherever possible. The graph thus created will be the same as replacing any crosses whose vertex lies on the face $(i,j)$ by a wrench with one vertex on $(i,j)$ and \emph{vice versa}. Also, the edge weights adjacent to $(i,j)$ are interchanged with their diametric opposites and the weight $x(i,j)$ is replaced by a certain binomial.

The graph thus produced is precisely the graph produced by the function $h'$ where $h'(i,j)=h(i,j)+2$ and $h'=h$ everywhere else. By induction, the matching polynomial of this graph is $f(n_0,i_0,j_0)$ for the initial conditions $h'$. The $f$ for the initial conditions $h$ is given by replacing $x(i,j)$ by the same binomial as before. Thus, $f(n_0,i_0,j_0)$ is precisely given by the matching polynomial of the graph for $h'$ with $x(i,j)$ replaced by the said binomial, which by the Urban Renewal Theorem is precisely the matching polynomial of the graph for $h$. \qedsymbol
 
\section{Proof II: Condensation} \label{Kuo}

In this section, we give a second proof of the main theorem, similar to the proof of \cite{Kuo}. In order to avoid lengthy tedious verifications, we will perform this proof with the face variables set equal to 1 throughout this section.

\subsection{The Condensation Theorem}

We first give some graph theoretic notations.

Let $G$ be a graph and $S \subseteq V(G)$. Let $\partial(S)$ denote the elements of $S$ that border vertices of $G$ not in $S$. If $G$ is bipartite, colored black and white, let $\delta(S)$ denote the number of black vertices of $S$ minus the number of white vertices. While $\delta(S)$ is only defined up to sign, we adopt the implicit assumption that, if $S_1$ and $S_2$ are both subsets of $V(G)$, then $\delta(S_1)$ and $\delta(S_2)$ are computed from the same coloring of $G$. Let $g(S)$ be the subgraph of $G$ induced by $S$. We abuse notation by writing $m(S)$ to mean $m(g(S))$. (Recall that all face variables have been set equal to 1, so $m(S)$ is a polynomial in the edge variables of $g(S)$.)

The following theorem is essentially due to Kuo and proven in \cite{Kuo}. (Kuo's paper proves this relation for many specific families of graphs but never states the general theorem. The below is my attempt to assimilate all of Kuo's cases into a general framework which will also encompass the results of this paper.)

\begin{Kuo}
Let $G$ be a bipartite planar graph and let the vertices of $G$ be partitioned into nine sets
$$V(G)=\C \sqcup \N \sqcup \NE \sqcup \E \sqcup \SE \sqcup \S \sqcup \SW \sqcup \W \sqcup \NW.$$
Assume that only vertices in the sets joined in figure~\ref{condense} border each other. (Vertices may also border other vertices in the same set.) 
Assume further that 
$$\delta(\NE)=\delta(\SW)=1,\ \delta(\SE)=\delta(\NW)=-1$$
and
$$\delta(\C)=\delta(\N)=\delta(\E)=\delta(\S)=\delta(\W)=0$$
Finally, assume that $\partial(\NE)$ and $\partial(\SW)$ are entirely black and $\partial(\NW)$ and $\partial(\SE)$ are entirely white.

\begin{figure}
\centerline{\scalebox{0.7}{\includegraphics{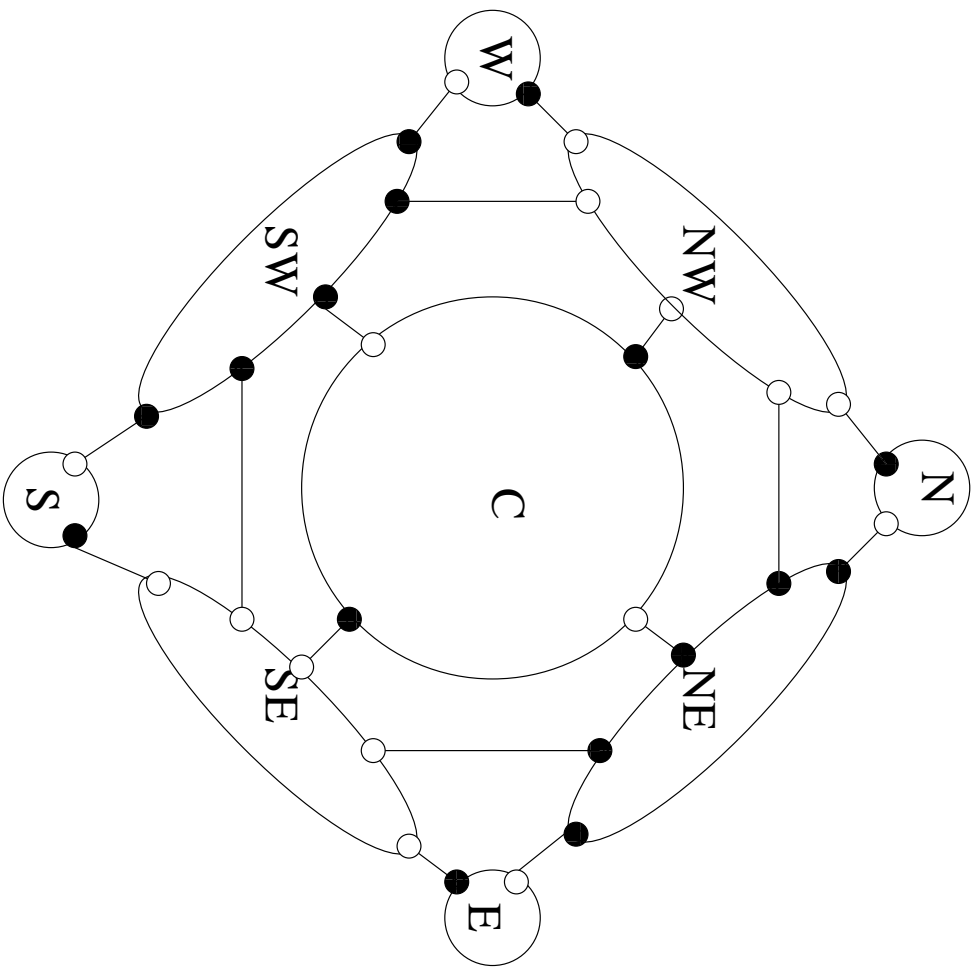}}}
\caption{The Connectivity and Coloring Hypotheses for the Condensation Theorem}\label{condense}
\end{figure}

Then we have
\begin{multline*}
m(G) m(\C)=m(\N \cup \NE \cup \NW \cup \C) m(\S \cup \SE \cup \SW \cup \C) m(\E) m(\W) + \\
m(\E \cup \NE \cup \SE \cup \C) m(\W \cup \NW \cup \SW \cup \C) m(\N) m(\S).
\end{multline*}
\end{Kuo}

We write $X$ to denote the collection of sets $\{ \N, \NE, \E, \SE, \S, \SW, \W, \NW \}$.

\begin{proof}
Let $G$ always denote a graph obeying the above hypotheses. We define a \emph{northern join} to be a set $\gamma$ of edges of $G$ of any of the following types
\begin{enumerate}
\item A path with one endpoint in $\NE$, the other endpoint in $\NW$ and the intermediate vertices in $\C$.
\item A single edge with one endpoint in $\NE$ and the other in $\NW$.
\item A pair of edges, one joining $\NE$ to $\N$ and the other joining $\NW$ to $\N$.
\end{enumerate}
We define an eastern, southern or western join analogously.

\begin{Lemma}
Let $M$ be a multi-set of edges of $G$ such that every member of $\C$ lies on two edges of $M$ and every other vertex lies on a single edge of $M$. Then we can write $M$ uniquely as a vertex disjoint union
$$M=\gamma_1 \sqcup \gamma_2 \sqcup M_{\C} \sqcup \bigsqcup_{\q \in X} M_{\q}$$
where
\begin{enumerate}
\item Either $\gamma_1$ is a northern join and $\gamma_2$ is a southern join, or $\gamma_1$ is a eastern join and $\gamma_2$ is a western join.
\item $M_{\C}$ is a disjoint union of cycles, entirely contained in $\C$. (We count a doubled edge as a cycle of length 2.)
\item $M_{\q}$ is a vertex disjoint set of edges contained entirely within $\q$.
\end{enumerate}
\end{Lemma}

\begin{proof}
Since every vertex is adjacent to either one or two elements of $M$, $M$ can be written uniquely as a disjoint union of cycles and paths. Moreover, the vertices of $\C$ are on cycles or are inner vertices of paths and the vertices of $V(G) \setminus \C$ are the ends of paths.

Let $M_{\C}$ be the cycles of $M$, $M_{\C}$ is vertex disjoint from the other edges of $M$. After deleting $M_{\C}$, the remaining graph is a disjoint union of paths, each of whose endpoints lie in one of the $\q \in X$ and all of whose interior vertices must lie in $\C$. Let $M_{\q}$ be the isolated edges connecting one point of $\q$ to another. After deleting all of the $M_{\q}$, what remains are paths as before which either have at least one interior point or whose endpoints lie in different $\q \in X$. Call the set of remaining edges $\Delta$.

Let $\q \in X$. We claim $\Delta \cap \q \subseteq \partial(\q)$. This is because, if $v \in \q \setminus \partial(\q)$, then there is a single edge $vw$ of $M$ containing $v$. We have $w \in \q$ by the definition of $\partial(\q)$. Then $vw$ is also the only edge containing $w$ and thus $vw$ is an isolated edge of $M$ lying entirely in $\q$. So $v \in M_{\q}$ and $\q \not\in \Delta$.

We know $M_{\NE}$ contains equally many black as white vertices as it is a disjoint union of edges. Thus $\Delta \cap \NE$ has one more black vertex than white, as we assumed $\delta(\NE)=1$. However, we just showed $\Delta \cap \NE \subseteq \partial(\NE)$, and $\partial(\NE)$ is entirely black. So there is a unique vertex $v_{\NE}$ in $\Delta \cap \NE$. Similarly, there are unique vertices $v_{\NW}$, $v_{\SE}$ and $v_{\SW}$, with $v_{\NE}$ and $v_{\SW}$ black and $v_{\NW}$ and $v_{\SE}$ white. Also by the same logic, $\#(\Delta \cap \q)$ is even for $\q \in \{\C, \NE, \NW, \SE, \SW \}$.

Now, the $v_{\q}$ must be endpoints of paths. The \emph{a priori} possible paths are the following, and their rotations and reflections:
\begin{enumerate}
\item A path from $v_{\NE}$ to $v_{\NW}$ through $\C$. (There are four possibilities in this symmetry class.
\item A path from $v_{\NE}$ to $v_{\SW}$ through $\C$. (There are two possibilities in this symmetry class.)
\item A single edge from $v_{\NE}$ to $v_{\NW}$. (There are four possibilities in this symmetry class.)
\item A single edge from $v_{\NE}$ to $N$. (There are eight possibilities in this symmetry class.)
\end{enumerate}

As $\#(\Delta \cap \N)$ is even, there must either be 0 or 2 paths of type (4) ending in $\N$, and similarly for $\E$, $\S$ and $\W$. This means we can not have a path joining $v_{\NE}$ to $v_{\SW}$. If such a path existed, $v_{\NW}$ could not be joined to $N$ or $W$, as that would create one path ending in that set. Also, $v_{\NW}$ could not be joined to $v_{\SE}$, as the paths $v_{\NE} v_{\SW}$ and $v_{\NW} v_{\SE}$ would cross. Similarly, we may not have a path joining $v_{\NW}$ to $v_{\SE}$. So the paths of type 2 do not occur.

It is now easy to see that $\Delta$ must decompose either as the union of a northern and a southern join or as the union of an eastern and a western join.
\end{proof}

Note that, if $\gamma_i$ ($i=1$ or $2$) passes through $\C$, it will contain an odd number of edges, as its endpoints are of opposite colors.

We now prove the result. Each side of the equation is a sum of products of edge variables, and each product corresponds to an $M$ meeting the conditions of the lemma; we will abuse notation and call this product $M$ as well. We will show each $M$ occurs with the same coefficient on both sides.

Specifically, let $k$ denote the number of cycles in $M_{\C}$ of length greater than 2. We claim that $M$ appears on each side of the equation with coefficient $2^k$. Even more specifically, we claim that $m(G) m(\C)$ contains $M$ with coefficient $2^k$. We also claim that, if $\gamma_1$ is a northern join and $\gamma_2$ a southern join, then $m(\N \cup \NE \cup \NW \cup \C) m(\S \cup \SE \cup \SW \cup \C) m(\E) m(\W)$ contains $M$ with a coefficient of $2^k$ and $m(\E \cup \NE \cup \SE \cup \C) m(\W \cup \NW \cup \SW \cup \C) m(\N) m(\S)$ does not contain $M$. If $\gamma_1$ is an eastern and $\gamma_2$ a western join then we claim the reverse holds.

Consider first the task of determining how many times $M$ appears in $m(G) m(\C)$. This is the same as the number of ways to decompose $M$ as $M=M(G) \sqcup M(\C)$ with $M(G)$ and $M(\C)$ matchings of $G$ and $\C$ respectively. 

Clearly, the edges with endpoints outside $\C$ must lie in $M(G)$. This means all the edges of $M_{\q}$, $\q \in X$, the edges of $\gamma_i$ if $\gamma_i$ does not pass through $\C$ and the final edges of $\gamma_i$ if $\gamma_i$ does pass through $\C$. This forces the allocation of all the edges of $\gamma_i$ because, when two edges of $M$ share a vertex, one edge must go in $M(G)$ and one in $M(\C)$. This forcing is consistent because $\gamma_i$ has an odd number of edges and the end edges must both lie in $M(G)$. Finally, we must allocate the edges of $M_{\C}$. All the cycles of $M_{\C}$ are of even length as $G$ is bipartite. In a cycle of length greater than 2, we may arbitrarily choose which half of its edges came from $M(G)$ and which from $M(C)$. Thus, we make $2^k$ choices.

Next, consider the coefficient of $M$ in $m(\N \cup \NE \cup \NW \cup \C) m(\S \cup \SE \cup \SW \cup \C) m(\E) m(\W)$. We must similarly write
$$M=M(\N \cup \NE \cup \NW \cup \C) \sqcup M(\S \cup \SE \cup \SW \cup \C) \sqcup M(\E) \sqcup M(\W).$$
We claim that if $\gamma_1$ is an eastern join and $\gamma_2$ a western, this coefficient is 0. There are three cases.

\setcounter{Case}{0}

\begin{Case}
$\gamma_1$ is a path joining $\NE$ to $\SE$ through $\C$. 
\end{Case}

Then every edge of $\gamma_1$ must lie in $M(\N \cup \NE \cup \NW \cup \C)$ or in $M(\S \cup \SE \cup \SW \cup \C)$ and which of these two it lies in must alternate. There are an odd number of edges in $\gamma_1$, so the final edges must lie in the same one of these two sets. But any edge with an endpoint in $\NE$ must lie in $M(\N \cup \NE \cup \NW \cup \C)$ and any edge with an endpoint in $\SE$ must lie in $M(\S \cup \SE \cup \SW \cup \C) \sqcup M(\E)$. \cont

\begin{Case}
$\gamma_1$ is a single edge with endpoints in $\NE$ and $\SE$. 
\end{Case}

But none of the $M(\q)$'s can contain such an edge. \cont

\begin{Case}
$\gamma_1$ is a union of two disjoint edges, one connecting $\NE$ to $\E$ and one connecting $\SE$ to $\E$.
\end{Case} 

But neither of these edge types can lie in any of the four $M(\q)$'s. \cont

Finally, we show that if $\gamma_1$ is a northern join and $\gamma_2$ a southern join that $M$ has coefficient $2^k$ in $m(\N \cup \NE \cup \NW \cup \C) m(\S \cup \SE \cup \SW \cup \C) m(\E) m(\W)$. As before, the edges of $M_{\q}$, $\q \in X$, are uniquely allocated to one of $M(\N \cup \NE \cup \NW \cup \C)$, $M(\S \cup \SE \cup \SW \cup \C)$, $M(\E)$ and $M(\W)$. If $\gamma_i$ does not pass through $\C$, its edges are also immediately forced. If $\gamma_i$ does pass through $\C$, its final edges are forced, thus forcing the others and again we have consistency as the two final edges are in the same set and $\gamma_i$ has an odd number of edges. Finally, each cycle of $M_{\C}$ of length greater than 2 can be allocated in two ways.

The case of $m(\E \cup \NE \cup \SE \cup \C) m(\W \cup \NW \cup \SW \cup \C) m(\N) m(\S)$ is exactly analogous. 
\end{proof}

\subsection{The Decomposition of $G$}

Let $G$ be a graph arising from the method of crosses and wrenches for some height function $h$ and some $(n_0, i_0, j_0)$. The purpose of this section is to describe a decomposition of the vertices of $G$ into nine disjoint sets. In the next section, we will show these sets obey the hypotheses of the previous section. For simplicity, we assume that $(n_0-2,i_0,j_0) \not\in \II$. 

Note that we have 
\begin{eqnarray*}
V(G_{(n_0-2,i_0,j_0)}) &=& V(G_{(n_0-1,i_0+1,j_0))}) \cap V(G_{(n_0-1,i_0-1,j))}) \\
                 &=& V(G_{(n_0-1,i_0,j_0-1)}) \cap V(G_{(n_0-1,i_0,j_0-1)})
\end{eqnarray*}
and
\begin{eqnarray*}
V(G_{(n_0,i_0,j_0)}) &=& V(G_{(n_0-1,i_0+1,j_0)}) \cup V(G_{(n_0-1,i_0-1,j)}) \\& & \cup V(G_{(n_0-1,i_0,j_0-1)}) \cup V(G_{(n_0-1,i_0,j_0-1)}).
\end{eqnarray*}
(Recall that $V(G)$ denotes the vertices of $G$.) 

\begin{figure}
\centerline{\scalebox{0.6}{\includegraphics{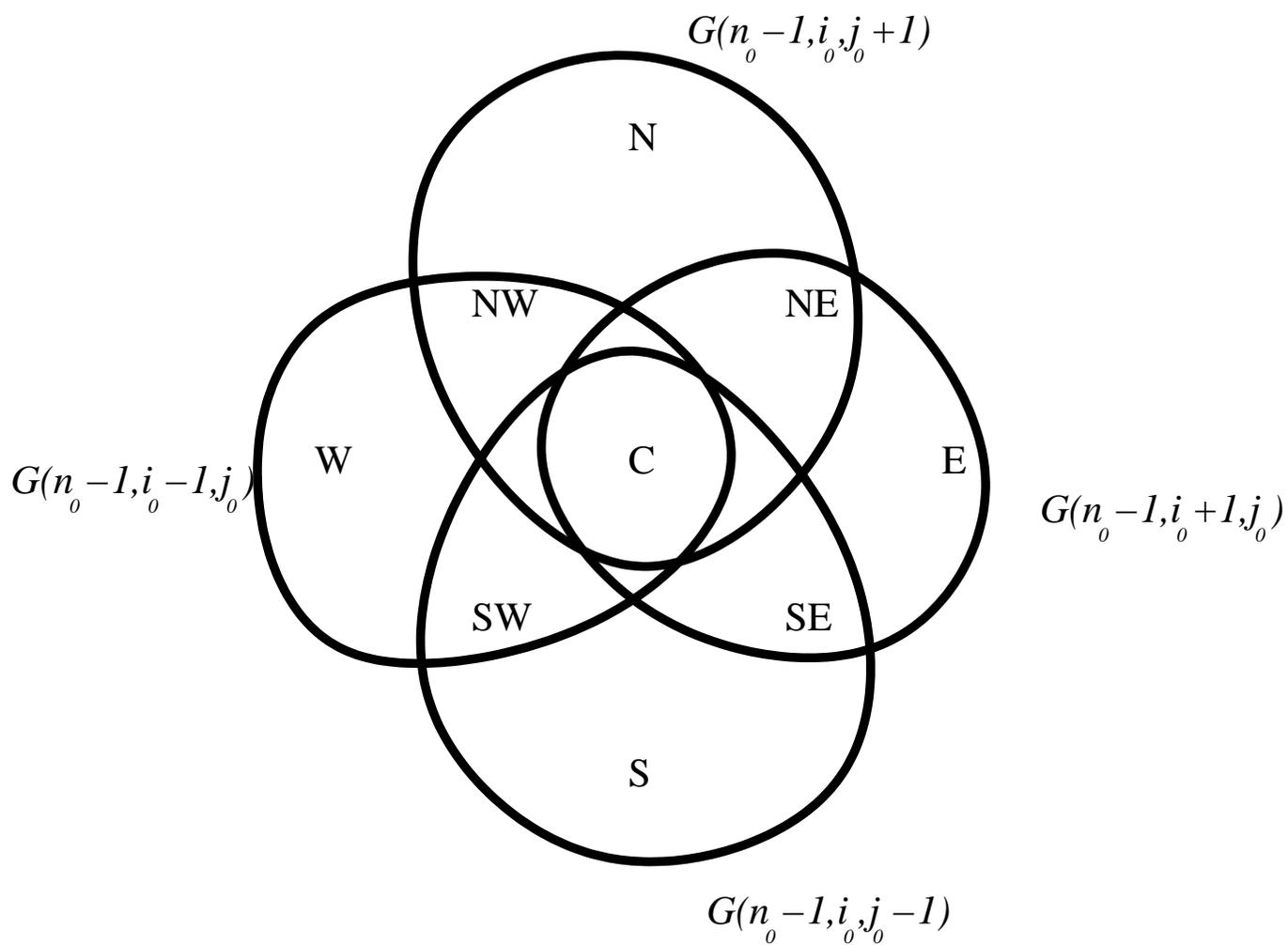}}}
\caption{Partitioning the Vertices of $G(n_0,i_0,j_0)$ into $9$ Sets} \label{partition}
\end{figure}

Thus, the four graphs $G(n_0-1,i_0 \pm 1, j_0)$ and $G(n_0-1, i_0,j_0 \pm 1)$ intersect as shown in the Venn diagram in Figure \ref{partition}, where the unlabelled cells indicate empty intersections. We decompose $V(G(n_0,i_0,j_0))$ into the nine sets $\C$, $\N$, $\NE$, $\E$, $\SE$, $\S$, $\SW$, $\W$ and $\NW$ as shown in the figure.

We now verify that the sets of vertices $\C$, $\N$, $\E$, $\W$, $\S$, $\NE$, $\NW$, $\SE$ and $\SW$ obey the hypotheses of the Condensation Theorem. We also show that $m(\E)=a(i_0+n_0-1,j_0)$, $m(\W)=c(i_0-n_0+1,j_0)$, $m(\N)=b(i_0,j_0+n_0-1)$, $m(\S)=d(i_0,j_0-n_0+1)$. These clearly establish the Main Theorem.

The following lemma is clear.

\begin{Lemma}\label{Ineq}
Abbreviate $p_{(n_0,i_0,j_0)}(i,j)$ by $p(i,j)$. Then 
\begin{align*}
(i, j) \in F_c(G_{(n_0-2,i_0,j_0)}) &\ iff\ & h(i,j)<p(i,j)-2 & & \\
(i, j) \in F_c(G_{(n_0-1,i_0+1,j_0)}) \setminus F_c(G_{(n_0-2,i_0,j_0)}) &\ iff\ & h(i,j)=p(i,j)-2 &\textrm{\ and\ }& i>i_0   \\
(i, j) \in F_c(G_{(n_0-1,i_0,j_0+1)} \setminus F_c(G_{(n_0-2,i_0,j_0)})) &\ iff\ & h(i,j)=p(i,j)-2 &\textrm{\ and\ }& j>j_0  \\
(i, j) \in F_c(G_{(n_0-1,i_0-1,j_0)}) \setminus F_c(G_{(n_0-2,i_0,j_0)}) &\ iff\ & h(i,j)=p(i,j)-2 &\textrm{\ and\ }& i<i_0  \\
(i, j) \in F_c(G_{(n_0-1,i_0,j_0-1)}) \setminus F_c(G_{(n_0-2,i_0,j_0)}) &\ iff\ & h(i,j)=p(i,j)-2 &\textrm{\ and\ }& j<j_0 
\end{align*}
(We are using our simplifying assumption that $(i_0,j_0) \in F_c(G_{(n_0-2,i_0,j_0)})$ to avoid a messy statement for this face.)
\end{Lemma} 

\begin{Proposition} \label{GeomN}
We have $m(\E)=a(i_0+n_0-1,j_0)$, $m(\W)=c(i_0-n_0+1,j_0)$, $m(\N)=b(i_0,j_0+n_0-1)$ and $m(\S)=d(i_0,j_0-n_0+1)$ and $\delta(\E)=\delta(\W)=\delta(\N)=\delta(\S)=0$. Also, all of the vertices of $\E$ that border vertices of $\NE$ are the same color. Calling this color white, the vertices of $\N$ bordering $\NE$, those of $\W$ bordering $\SW$ and those of $\S$ bordering $\SW$ are white. The vertices of $\E$ bordering $\SE$, those of $\S$ bordering $\SE$, those of $\N$ bordering $\SW$ and those of $\W$ bordering $\SW$ are black.
\end{Proposition}

\begin{proof}
We do the case of $\E$, the others are similar. Let $v \in \EE$ and let $(i,j)$ be a closed face of $G_{(n_0-1,i_0+1,j_0)}$ bordering $v$ (by the definition of $\E$, such a face exists); by Lemma~\ref{Ineq}, $h(i,j)=p(i,j)-2$. Now, $(i,j)$ is not a closed face of $G_{(n_0-1,i_0,j_0+1)}$. Thus, by Lemma~\ref{Ineq}, $j\leq j_0$. But, similarly, $j \geq j_0$ so $j=j_0$. Moreover, the faces that $v$ borders which are not of the form $(i,j_0)$ must not lie in $G$ at all, so we have $h(i,j \pm 1)=h(i,j)+1$. 

Putting this data into the Crosses and Wrenches algorithm, we have that $\E$ looks like figure~\ref{ELayout}, where the number of hexagons and the length of the dangling path can vary. The highlighted vertices are those that border $\NE$ and $\SE$.

\begin{figure}
\centerline{\includegraphics{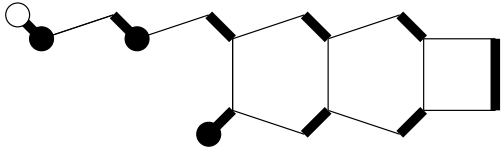}}
\caption{The Geometry and Coloring of $\E$}\label{ELayout}
\end{figure}

Clearly, $\E$ has a unique matching $M$, shown in bold. Let $i$ be the largest value for which $(i,j_0)$ is a closed face of $G$, the only weighted edge in $M$ is the edge separating $(i,j_0)$ and $(i,j_0+1)$. As $h(i,j_0)=(n_0-1)+(i_0+1-i)-(j_0-j_0)=n_0+i_0-i$ and $h(i+1,j_0)=h(i,j_0)+1=n_0-i_0+i+1$ (as otherwise $(i+1,j_0)$ we be a closed face) we see that this is the edge $a(i_0+n_0-1,j_0)$. As $\E$ has a matching, $\delta(\E)=0$.

All that is left to show is that the vertices of $\N$ that border $\NE$ and the vertices of $\E$ that do are the same color. It is equivalent to show that the vertex of $b(i,j+n-1)$ and that of $d(i+n-1,j)$ closer to $\NE$ are the same color. But the path joining them is the one we proved in Proposition \ref{Boundary} had an odd number of edges.
\end{proof}

\begin{Lemma}
Let $v \in \NE$. Then exactly one of the following is true.
\begin{enumerate}
\item $v$ lies in a wrench oriented NE-SW (i.e., this orientation: $\wrench$), the other end of which is also in $\NE$. \\
or
\item $v$ lies on the outer perimeter of $G_{(n_0,i_0,j_0)}$. 
\end{enumerate}
\end{Lemma}

\begin{proof} 
Let $(i,j)$, $(i+1,j)$, $(i,j+1)$, $(i+1,j+1)$ be the faces surrounding $v$ ($v$ may only be adjacent to three of them) and let $I$ be the set of these four faces.. By Lemma~\ref{Ineq}, either $i$ or $i+1 > i_0$ so $i \geq i_0$. Similarly, $j \geq j_0$. There are several cases.

First, suppose $h(i,j)=p(i,j)-2$. Then $h(i',j') \geq p(i',j')-2$ for $(i',j') \in I$. If $h(i',j')=p(i',j')-2$ for all $(i',j') \in I$, then $v$ lies on a NE-SW wrench both ends of which are in $\NE$. If some $h(i',j')>p(i',j')-2$ then $(i',j')$ is not a closed face of $G_{(n_0,i_0,j_0)}$ and one can check that in every case $v$ is adjacent to $(i',j')$, so $v$ is on the outer perimeter.

Now, suppose that $h(i,j)<p(i,j)-2$. Then $(i,j) \in F_c(G_{(n_0-2,i_0,j_0)})$. The only way that $v$ could be in $\E$ then is if $v$ were at the NE end of a NE-SW wrench so that it didn't border $(i,j)$. Then $(i+1,j)$ does border $v$, so $h(i+1,j) \geq p(i+1,j)-2$. The only way this is possible is if $h(i,j)=p(i,j)-4$, $h(i+1,j)=h(i,j+1)=p(i,j)-3$ and $h(i+1,j+1)=p(i,j)-2=p(i+1,j+1)$. Then $(i+1,j+1) \not \in F_c(G_{(n_0,i_0,j_0)})$ and $v$ lies on $(i+1,j+1)$, so again $v$ is on the perimeter.

Finally, suppose that $h(i,j)>p(i,j)-2$. Then $h(i',j') > p(i',j')-2$ for all $(i',j') \in I$. But then none of the $(i', j')$ are closed faces of $G_{(n_0,i_0,j_0)}$ and $v$ is not a vertex of $G_{(n_0,i_0,j_0)}$, a contradiction.
\end{proof}

\begin{Proposition} $\NE$ and $\SW$ have one more black vertex than white, and \emph{vice versa} for $\SE$ and $\NW$.
\end{Proposition}

\begin{proof}
We discuss the case of $\NE$, the others are similar.

Pair off all of the vertices of $\E$ for which the first case of the previous proposition holds with the other end of their wrench. This removes equal numbers of black and white vertices. What remains is a path, which we showed in Proposition~\ref{GeomN} has each end black.
\end{proof}

\begin{Proposition} 
$\partial(\NE)$ and $\partial(\SW)$ are entirely black and \emph{vice versa} for $\partial(\SE)$ and $\partial(\NW)$. Moreover, only the adjacencies permitted by the hypotheses of the Condensation Theorem can occur.
\end{Proposition}

\begin{proof}
Besides the endpoints of the path in the previous proposition, $\partial(\NE)$ must be made up entirely of vertices from the wrenches in the previous proposition. However, in fact, only the south-west vertex of a wrench will be able to lie in $\partial(\NE)$, and these are all the same color.

This deals with the coloration of the boundaries. This immediately means that $\NE$ can not border $\SW$. We saw in Lemma~\ref{GeomN} that $\NE$ only borders $\N$ and $\E$. These and the symmetrically equivalent restrictions gives the required restrictions on what may border what.
\end{proof}

\begin{Proposition} 
$\C$ has equally many black and white vertices.
\end{Proposition}

\begin{proof}
$\C=G_{(n_0-2,i_0,j_0)}$. By induction on $n_0$, we can assume that we already know that $f(n_0-2,i_0,j_0)=m(\C)$. In particular, $\C$ has at least one perfect matching.
\end{proof}

We have now checked all of the conditions to apply the Kuo Condensation Theorem.

\subsection{A Consequence of the Condensation Proof}

The Condensation proof of the main theorem, although it involves more special cases, seems to have at least one consequence which does not follow from the Urban Renewal proof.

\begin{Proposition}
In the relation
\begin{eqnarray*}
f(n,i,j) f(n-2,i,j) &=& a(i+n-1,j)c(i-n+1,j)f(n-1,i,j+1)f(n-1,i,j-1) +   \\
           & &  b(i,j+n-1)d(i,j-n+1)f(n-1,i+1,j)f(n-1,i-1,j)
\end{eqnarray*}
each monomial appearing on the right hand side appears in exactly one of the two terms on the right, and appears with the coefficient $2^k$ for some $k$.
\end{Proposition}

\begin{proof}
This just says that every combination of edge $M$ as in the proof of the Kuo Condensation Theorem occurs either only in 
$$m(\N \cup \NE \cup \NW \cup \C) m(\S \cup \SE \cup \SW \cup C) m(\E) m(\W)$$
or only in
$$m(\E \cup \NE \cup \SE \cup \C) m(\W \cup \NW \cup \SW \cup C) m(\N) m(\S)$$
and appears with coefficient $2^k$ for some $k$; this follow from the proof of the Kuo Condensation Theorem.
\end{proof}

\section{Some Final Observations} \label{Final}

\subsection{Setting Face or Edge Variables to One}

\begin{Proposition}
If we replace all of the face variables in $f(n,i,j)$ by 1, every coefficient in $f(n,i,j)$ will still be 1.
\end{Proposition}

\begin{proof}
We are being asked to show that we can recover a matching of $G_{n,i,j}$ (henceforth, $G$) from knowing which weighted edges were use in the matching. Since the unweighted edges of $G$ are vertex disjoint, knowing which weighted edges are used determines all the edges used, and hence the matching.
\end{proof}

The corresponding version for faces is true but more difficult, and its proof has more interesting consequences.

\begin{Proposition}\label{facedetermine}
If we replace all of the edge variables in $f(n,i,j)$ by 1, every coefficient in $f(n,i,j)$ will still be 1.
\end{Proposition}

\begin{proof}
We are now being asked to show that we can recover a matching $M$ of $G$ from just knowing the face exponents. By the previous result, it is enough to show that we can determine whether or not a weighted edge appears in $M$. Recall the notation of Section~\ref{gwofs}: for $f$ a face of $\GG$, $\epsilon(f)$ is the exponent of $x_f$ and, for $e$ an edge of $\GG$, $\delta(e)$ is the exponent of $e$.

We give a formula for $\delta_e$ in terms of certain $\epsilon(f)$'s in the case where $e$ is of the form $(a,i,j)$, similar formulas for $\delta(b,i,j)$, $\delta(c,i,j)$ and $\delta(d,i,j)$ can be found and will be stated later in the section. My thanks to Gabriel Carrol for suggesting this simple formulation and proof of the lemma (in a somewhat different context.)
\begin{Lemma}
Assume that $G_{n,i,j}$ contains the edge $(a,i_0,j_0)$. In any matching of $G$, we have
$$\delta(a,i_0,j_0)=-\sum_{\substack{(n,i,j) \in \II \\ n+i+j<i_0+j_0+1 \\n+i-j < i_0-j_0+1}} \epsilon(n,i,j).$$
\end{Lemma}

\begin{proof}
Let $t$ denote a formal variable independent of all other variables.

Let $f$ be defined by 
\begin{eqnarray*}
f(n,i,j) f(n-2,i,j) &=& a(i+n-1,j)c(i-n+1,j)f(n-1,i,j+1)f(n-1,i,j-1) +   \\
           & &  b(i,j+n-1)d(i,j-n+1)f(n-1,i+1,j)f(n-1,i-1,j)
\end{eqnarray*}
for the initial conditions $x(i,j)$ on $\II$. 

Define $\tilde{f}$ by
\begin{eqnarray*}
\tilde{f}(n,i,j) &=& t f(n,i,j) \qquad n+i+j<i_0+j_0+1\ \textrm{and}\ n+i-j<i_0-j_0+1 \\
\tilde{f}(n,i,j) &=& f(n,i,j) \qquad \textrm{otherwise.}
\end{eqnarray*}
Define $\tilde{a}(i_0,j_0)=ta(i,j)$ and $\tilde{a}(i,j)=a(i,j)$ otherwise. Then we claim that
\begin{eqnarray*}
\tilde{f}(n,i,j) \tilde{f}(n-2,i,j) &=& \tilde{a}(i+n-1,j)c(i-n+1,j)\tilde{f}(n-1,i,j+1)\tilde{f}(n-1,i,j-1) +   \\
           & &  b(i,j+n-1)d(i,j-n+1)\tilde{f}(n-1,i+1,j)\tilde{f}(n-1,i-1,j).
\end{eqnarray*}

This may be checked by checking all nine cases for how $n+i+j$ and $n+i-j$ compare to $i_0+j_0+1$ and $i_0-j_0+1$, and noting that $n+i+j=i_0+j_0+1$ and $n+i-j=i_0-j_0+1$ imply that $i+n-1=i_0$ and $j=j_0$.

Now, let us consider the coefficient of $t$ in any term of $\tilde{f}(n,i,j)$. Since $(a,i_0,j_0)$ appears in $G_{n,i,j}$, we must have $n+i+j \geq i_0+j_0+1$ and $n+i-j \geq i_0-j_0+1$. Thus, $\tilde{f}(n,i,j)=f(n,i,j)$ and no $t$ appears. On the other hand, $\tilde{f}$ can also be obtained from $f$ by substituting $t x(i,j)$ for $x(i,j)$ for every $(n,i,j) \in \II$ such that $n+i+j < i_0+j_0+1$ and $n+i-j < i_0-j_0+1$ and substituting $t a(i_0,j_0)$ for $a (i_0,j_0)$. So we see that the exponent of $t$ is
$$\delta(a,i,j)+\sum_{\substack{ (n,i,j) \in \II \\ n+i+j<i_0+j_0+1 \\n+i-j < i_0-j_0+1}} \epsilon(n,i,j).$$

So this sum is zero and we are done.
\end{proof}

So the face exponents determine the edge exponents and (by the previous Theorem) determine the matching.
\end{proof}

One can find three more formulas for $\delta(a,i,j)$ and four corresponding formulas for each of $\delta(b,i,j)$, $\delta(c,i,j)$ and $\delta(d,i,j)$. Specifically, let 
\begin{eqnarray*}
P &=& \{ (n,i,j) \in \II : n+i+j < i_0+j_0+1,\ n+i-j < i_0-j_0+1 \} \\
Q &=& \{ (n,i,j) \in \II : n+i+j < i_0+j_0+1,\ n+i-j \geq i_0-j_0+1 \} \\
R &=& \{ (n,i,j) \in \II : n+i+j \geq i_0+j_0+1,\ n+i-j \geq i_0-j_0+1 \} \\
S &=& \{ (n,i,j) \in \II : n+i+j \geq i_0+j_0+1,\ n+i-j < i_0-j_0+1 \} \\
\end{eqnarray*}

Then
$$\delta(a,i,j) = - \sum_{(n,i,j) \in P} \epsilon(n,i,j) 
              =   \sum_{(n,i,j) \in Q} \epsilon(n,i,j) 
              = 1-  \sum_{(n,i,j) \in R} \epsilon(n,i,j) 
              =   \sum_{(n,i,j) \in S} \epsilon(n,i,j) $$

The same result holds for
\begin{enumerate} 
\item $\delta(b,i,j)$, if we define $P$, $Q$, $R$ and $S$ by comparing $n \pm i +j+1$ to $\pm i_0 +j_0$.
\item $\delta(c,i,j)$, if we define $P$, $Q$, $R$ and $S$ by comparing $n - i \pm j+1$ to $-i_0 \pm j_0$.
\item $\delta(d,i,j)$, if we define $P$, $Q$, $R$ and $S$ by comparing $n \pm i - j+1$ to $\pm i_0 - j_0$.
\end{enumerate}

A geometrical description can be given for these rules. Let $e$ be a weighted edge of $G$. Draw four paths in $G$ by starting at $e$ and alternately turning left and right. (There are two ways to leave $e$ and two choices for which turn to make first.) These divide the faces of $G$ into four regions $P$, $Q$, $R$ and $S$. If the edge $e$ is present, then the sums $\sum_{(n,i,j) \in P} \epsilon(n,i,j)$,  $\sum_{(n,i,j) \in Q} \epsilon(n,i,j)$,   $\sum_{(n,i,j) \in R} \epsilon(n,i,j)$ and  $\sum_{(n,i,j) \in S} \epsilon(n,i,j)$ will yield a 1, two -1's and a 0. If $e$ is absent, one of these sums will be 1 and the other three will be 0.

Interestingly, if one follows the operation of the previous where $e$ is an unweighted edge, one gets precisely the same conclusions except that the interpretation of the results is changes: getting the sums $(1,0,0,0)$ (in some order) now means that $e$ is \emph{present} and $(0,1,-1,1)$ means $e$ is \emph{absent}.

\subsection{Height Functions for Crosses and Wrenches Graphs}

Let $G$ be a bipartite connected planar graph with a fixed labeling of its vertices as black and white and let $w$ be a function assigning a real number to each edge of $G$. We define a \emph{Propp height} to be a function $H$ from the faces of $G$ to the real numbers with the following property: Let $f$ and $f'$ be faces of $G$ separated by an edge $e$ such that, when one stands on $f$ and looks toward $f'$, the white vertex of $e$ is on the left. Then $H(f)-H(f')=w(e)+(0 \textrm{\ or\ } -1)$. Define two Propp heights $H$ and $H'$ to be equivalent if $H(f)=H'(f)+c$ for some constant $c$ independent of $f$. (The term Propp height is my own, to distinguish them from the heights which occur through out this paper.)

In section 3 of \cite{Propp1}, Propp essentially proved
\begin{Heights}
For $G$ as above, one can find $w$ such that there is a bijection between equivalence classes of Propp heights for $w$ and matchings of $G$. The bijection is such that, if $H$ corresponds to $M$, the edge $e$ occurs in $M$ iff $H(f)=H(f')+(w(e)-1)$ where $e$ separates $f$ and $f'$.
\end{Heights}

These Propp heights have been very useful in studying statistical properties of random matchings, see \emph{e.g.} \cite{CEP}. For crosses and wrenches graphs, the Propp height has an extremely simple description.
\begin{Proposition}
Let $\GG$ be a crosses and wrenches graph. We may take $w$ in the previous theorem to be given by $w(e)=1/4$, if $e$ is a horizontal or vertical edge and $w(e)=1/2$, if $e$ is a diagonal edge.
\end{Proposition}

\begin{proof}
We just must verify that, for every vertex $v$ of $\GG$, $\sum_{e \ni v} w(e)=1$, where the sum is over $e$ all incident on $v$. This is clear.
\end{proof}

In this context, the condition at infinity for infinite completions can be stated as $H(i,j)=(|i|-|j|)/4$ for all but finitely many $(i,j)$.

\subsection{Generating Random Matchings}

Let $G$ be a crosses and wrenches graph. Suppose that we want to randomly generate a matching of $G$, with our sample drawn from all matchings of $G$ with uniform probability. Such sampling has produced intriguing results and conjectures in the cases of Aztec diamonds and fortresses. 

It is possible to do so in time $O(|G|^2)$, where $|G|$ can be any of $|F(G)|$, $|E(G)|$ or $|V(G)|$ as all of these only differ by a constant factor for $G$ a crosses and wrenches graph. More specifically, if $G=G_{(n_0,i_0,j_0)}$, one can do so in time $O(|\C_{(n_0,i_0,j_0)} \cap \UU|)$. (In most practical cases, this is closer to $|G|^{3/2}$.

We give only a quick sketch; the method is a simple adaptation of the methods of \cite{Propp2}. Let $G=G_{(n_0,i_0,j_0)}$ arise from a height function $h$ and let $F=F(G)$. Throughout our algorithm, $x$ will denote a function $F \to \RR$. Our algorithm takes as input $h$ and $x$ and returns a random matching, where the probability of a matching $M$ being returned is proportional to $\prod_{f \in F} x(f)^{\delta(f)}$. We describe our output as a list of the weighted edges used in $M$. At the beginning of the algorithm, we take $x(f)=1$ for all $f$.

\textbf{Step 1:} 
If $h(i_0,j_0)=n_0$, return $\{ \}$ and halt.

\textbf{Step 2:}
Find $(i,j) \in F$ such that $h(i,j)=h(i \pm 1, j \pm 1)-1$ for all four choices of the $\pm$. 

\textbf{Step 3:}
Set $x'(i,j)=(x(i+1,j) x(i-1,j)+x(i,j+1) x(i,j-1))/x(i,j)$ and $x'(f)=x(f)$ for all other $f \in F$. Set $h'(i,j)=h(i,j)+2$ and $h'(f)=h(f)$ for all other $f$.

\textbf{Step 4:}
Run the algorithm with input $h'$ and $x'$, let the output be $M'$.

\textbf{Step 5:}
For shorthand, set $a=(i+1+h(i,j),j,a)$, $b=(i,j+1+h(i,j),b)$, $c=(i-1-h(i,j),j,c)$ and $d=(i,j-1-h(i,j),d)$. If $M'$ contains one of $\{ a,b,c,d \}$, return M and halt. If $M'$ contains two of $\{ a,b,c,d \}$, return $M' \setminus \{ a,b,c,d \}$ and halt. If $M'$ contains none of $\{ a,b,c,d \}$, return $M \cup \{ a,c \}$ with probability $x(i,j+1) x(i,j-1)/x(i,j) x'(i,j)$ and return $M \cup \{ b,d \}$ with probability $x(i+1,j) x(i-1,j)/x(i,j) x'(i,j)$. All other cases are impossible.

In order to not do a time consuming search at step 2, one can keep a reverse look-up table which, given $h$, record all $(i,j)$ with $h(i,j)=h$ and is updated in step 3 when $h$ is replaced by $\tilde{h}$. Note that $x(f)$ will always be integer valued.

\end{document}